\documentstyle [12pt] {article}
\begin{document}
\newfont{\fraktur}{cmbx10}
\newcommand{\C}{{\bf C}}
\newcommand{\HH}{{\bf H}}
\newcommand{\K}{{\bf K}}
\newcommand{\N}{{\bf N}}
\newcommand{\R}{{\bf R}}
\newcommand{\Z}{{\bf Z}}
\newcounter{roman}
\newtheorem{definition}{Definition}[section]
\newtheorem{lemma}[definition]{Lemma}
\newtheorem{theorem}[definition]{Theorem}
\newtheorem{corollary}[definition]{Corollary}
\newtheorem{proposition}[definition]{Proposition}
\newtheorem{remark}[definition]{Remark}
\parindent=0pt
\hfuzz=3pt
\begin{titlepage}
\setcounter{page}{-1}
\title{Covariant SPDEs and Quantum Field Structures}
\author{C.Becker\\
Institut f\"ur Mathematik\\Ruhr-Universit\"at Bochum\\
44780 Bochum\\Germany \and
R.Gielerak\\
P.\L ugiewicz\\
Institute of Theoretical Physics\\University of Wroc\l aw\\
50-204 Wroc\l aw\\
Poland}
\date{June 26, 1996}
\maketitle

\begin{abstract}
\par\noindent
Covariant stochastic partial differential equations are studied in any
dimension.
A special class of such equations is selected and it is proven that
the solutions can be analytically continued to Minkowski space-time
yielding tempered Wightman distributions which are covariant, obey the
locality axiom and a weak form of the spectral axiom.
\end{abstract}

PACS numbers:
\thispagestyle{empty}
\end{titlepage}
\newpage

\tableofcontents
\thispagestyle{empty}

\newpage


\section{Introduction}
The connection between scalar generalized random fields which are Markov and
Euclidean invariant and scalar quantum fields played a crucial role in the
development of constructive quantum field theory \cite{Glimm-Jaffe,Simon}.
Symanzik \cite{Symanzik} first pointed out this connection for the free field
and Nelson \cite{Nelson1,Nelson2} developed some general machinery to
construct quantum fields from Euclidean invariant Markov fields.
Multi-component Gaussian generalized random fields which are Markov and
invariant under the Euclidean group might  play a role similar to that of
the free scalar field \cite{Ito,Wong1,Wong2,Yao}.
A simple example for such covariant random fields is given by
infinitely divisible random fields \cite{BergForst}.
It seems that these fields are too singular:
perturbations by local multiplicative
functionals as in the standard constructive quantum field theory
approach should lead to a very serious ultraviolet divergence problems;
nevertheless there is another
constructive approach  which was initiated in
\cite{Albphyslettb1,Albphyslettb2,Albnoise,Albcomm} and in the
following papers \cite{AlbevI,AlbevII,Osipov1,Osipov2}.
In all the above-mentioned papers it is essentially needed that a real
vector space of dimension $D=1,2,4,8$ can be given the structure
of a division algebra so that the Laplace operator
$\bigtriangleup _{D}=\sum_{i=1}^{D} \frac{\partial
^{2}}{\partial x^{2}_{i}}$ can be factorized as a product of two first-order
covariant elliptic differential operators $\partial$ and
$\overline{\partial}$. One can then consider an equation of the form
\begin{equation}
\partial A = \eta
\end{equation}
where $\eta$ is suitably chosen noise.
The solution of this equation, which can be computed explicitly, is again
a covariant Markovian generalized random field.
The moments of this generalized random field can be analytically continued to
Minkowski space-time, yielding covariant system of Wightman distributions which
obey the locality axiom and a weak form of the spectral axiom
\cite{Bogoliubov,Jost,Streater}.
Moreover, if the noise $\eta$ contains a nonzero Poisson
piece the corresponding system of Wightman functions is not quasi-free
(non Gaussian).
\par
In the present paper we shall consider an equation of the type
\begin{equation}
\label{fundamental_equation}
{\cal D}A=\eta
\end{equation}
in arbitrary space-time dimension $D \geq 2$ and where ${\cal D}$
is an arbitrary
covariant differential operator of any order.
\par
It is among the main objectives of the present paper
to demonstrate that the existence of division algebras in the particular
dimensions $D=1, 2, 4$, and $8$ is not essential and that in any dimension
the covariant
Markovian generalized random field $A$ can be constructed by solving
equation $(2)$ with suitable ${\cal D}$ and $\eta$. Morever it
will be shown that it is a generic
property of a large class of such equations that the moments of the
random field
$A$ can be analytically continued to Minkowski space giving a set of
tempered Wightman distributions which are covariant and which fulfill the
locality  axiom and a weak form of the spectral axiom.
\par
The essential problem behind these constructions is to decide
whether a reflection-positive non-Gaussian covariant generalized random field
$A$ can be obtained from equation $(2)$.
Unfortunately, the authors have obtained some negative results which will
be published in forthcoming papers.
One of the negative results, obtained by the second and the third author
\cite{Romek2}, directly relates to solutions of equation
(\ref{fundamental_equation}).
\par
Let ${\tau}$ be a real finite-dimensional representation of the
orthogonal group ${\bf O(D)}$ with  $D \geq 2$ and let
$\tau = \oplus _{\alpha}
\tau _{\alpha}$ be its decomposition into irreducible
representations. Let $A=(A_{\alpha})_{\alpha}$ denote
the corresponding decomposition of the $\tau$-covariant field
$A$ which solves equation (\ref{fundamental_equation}).
Then for
any $\alpha$ such that $\hbox{dim}\,\tau _{\alpha }>1$ the corresponding
Euclidean
field
$A_{\alpha}$ is not reflection-positive.
\par
However, this does not exclude the possibility that on
the scalar sectors of $\tau =\oplus_{\alpha} \tau_{\alpha}$
(i.e.\ on the subspaces where $\hbox{dim}\, \tau _{\alpha}=1$) or on a certain
subspaces of ${\cal S}(\R^{D})\otimes \R^{\hbox{dim}\, \tau _{\alpha}}$
reflection positivity holds. Moreover the possibility of passing
to complex representations of ${\bf O(D)}$ is not covered by this
negative result.
\par
Another negative result, obtained by the first author
\cite{Becker1,Becker2}, is that there is essentially no multi-component
generalized random field which is covariant with respect to some
representation
$\tau=\oplus_\alpha \tau_\alpha$ of {\bf O(D)} and which is
reflection-positive on the whole test function space unless at least one of the
$\tau_\alpha$ is trivial.
However, one can construct fields of ultra-local type which lead to a trivial
Hilbert space.
In order to obtain a field which is non-trivial from the point of view of
physics, one has to restrict the test function space which corresponds to
fixing some gauge.
\par
One of the main conclusions from these negative results is that
in order to obtain reflection-positive bosonic random fields of spin higher
than zero one has to demand covariance only with respect to {\bf SO(D)}
instead of {\bf O(D)}.
This is not in conflict with the Wightman
axioms \cite{Bogoliubov,Jost,Streater} since in the axiomatic framework
the existence of Euclidean fields is
not required at all and moreover the covariance of the corresponding set of
Schwinger functions is demanded with respect to ${\bf SO(D)}$ only.
For a more detailed exposition we refer to forthcoming papers
\cite{Becker1,Becker2,Romek2}.
For a construction of Euclidean fields of arbitrary spin in an axiomatic
framework we refer to \cite{Ozkaynak}.
\par
It seems to be an intrinsic property of gauge fields that the conditions of
positivity, covariance and locality are all together not compatible with
local gauge invariance \cite{Strocchi1,Strocchi2}.
In view of this, we expect that some of the
models produced by the methods described in the present paper,
though they are not reflection-positive,
could find applications in problems of quantum field theory
of gauge type with indefinite metrics. This is the second motivation for
the present and some forthcoming papers \cite{Romek1,Piotrek}.
\smallskip
\par
{\bf Organisation of the paper}
\par
Although the proper mathematical language for the material presented
in this paper is the language
of vector bundles over $\R^{D}$ and equivariant differential operators
of first order we decided to present our results in a more elementary way in
order to make them easily accessible to a wider audience.
In section 2 we fix
the notation and mention some elementary
results which some of the readers probably know.
The main result of the paper is contained in
section $3$:
Assume that
${\cal D}$
has admissible mass spectrum (see below for the definition) and that
$\eta $ is white noise that possesses
all moments.
Then there exist tempered covariant
distributions supported in the forward cone  such that their Laplace-Fourier
transforms are equal to the moments of $A$ regarded as functions of the
difference variables at positive time.
Finally, in the last section we
present some particular examples in three-dimensional space resulting from the
lowest-dimensional real representations of the group ${\bf SO(3)}$.
Models describing the
interaction between scalar fields
and vector fields that we call Higgs${}_{3}$-like
models and models describing two interacting vector fields are also
presented the last section.

\section{Random Fields as Solutions of Covariant SPDEs}
\subsection{Covariant First-Order Differential Operators}

An important concept in physics is the concept of covariance, i.e. the
fact that the form of an equation does not change under
suitable coordinate transformations.
There is a lot of literature on this subject, cf.\
\cite{Fonda,GelfandShapiro,Ljubarski,WaveEquations}.
In this section we shall investigate covariant first-order differential
operators acting on $C^{\infty}$-functions $\R^D\to \K^N$,
$\K \in \{\R,\C\}$.
We assume that a representation of a Lie group
${\bf G}\!\subseteq\! {\bf GL(D)}$
is acting on $\K^N$.
In our applications we shall mainly study the case ${\bf G}={\bf SO(D)}$,
which is motivated by our intention to produce covariant models in the
framework of Euclidean quantum field theory.
\smallskip
\par
Let us, first of all, collect some basic definitions and facts.
\begin{proposition}
\label{covoperators}
${ }$
\par
Let
$B_1,\ldots,B_D$ be matrices $\in{\cal M}_{N\times N}(\K)$
and put
$B=(B_1,\ldots,B_D)$.
Let $E\in{\cal M}_{N\times N}(\K)$ denote the unit matrix.
\par
We consider
the first-order operator
\begin{equation}
\label{covoperator}
{\cal D}_B=\sum_{j=1}^D\; B_j\; {\partial\over\partial x_j} + m\; E \ \ ,
\quad m\in\R
\end{equation}
acting on the space of $C^{\infty}$-functions
$\R^D\to \K^N$.
Let $T_g$ denote the action of the representation $\tau$ on functions
$f\in C^\infty(\R^D,\K^N)$:
\begin{equation}
T_g\, f(x) = \tau(g) f(g^{-1}x),\quad g\in{\bf G}\ \ \ .
\end{equation}
The following statements are equivalent:
\smallskip
\begin{itemize}
\item [(a)]
The form of ${\cal D}_B$ does not change if we make a coordinate transformation
in $\R^D$ :
$x\mapsto g x $, $g\in {\bf G}$,
and simultaneously a coordinate transformation in $\K^N$:
$y \mapsto \tau(g)\, y$.
\smallskip
\item [(b)]
${\cal D}_B$ commutes with $T_g$ :
\begin{equation}
\big[  {\cal D}_B , T_g \big] = 0 \quad \forall g \in {\bf G}\ \ \ .
\end{equation}
\smallskip
\item [(c)]
\begin{equation}
\label{defcovariant}
 \sum_{k=1}^D\; g_{jk}\; \tau(g)\; B_k\; \tau(g^{-1}) = B_j
 \qquad\forall j\in\{1,\ldots,D\}\qquad\forall g \in {\bf G}
\end{equation}
where $g_{jk}$ are the components of $g  \in {\bf G}$.
\end{itemize}
\end{proposition}
\medskip

Note, that instead of taking the operator $m\cdot E$ in $(3)$ we can
take any matrix $M$ belonging to the commutant of the representation
$\tau$.
\begin{definition}
If ${\cal D}_B$ fulfills one (and hence all) of the conditions in proposition
\ref{covoperators}, it will be called covariant with respect to the
representation $\tau$.
\par
The set of all operators that are covariant with respect to $\tau$ will
be denoted by ${\rm Cov}(\K^N,\tau)$.
\end{definition}
\smallskip
\par
Note that if $\tau(g) \in {\bf O(N)} \quad \forall g \in {\bf G}$ and
if
$(B_1,\ldots,B_D)$ defines a covariant operator with respect to $\tau$
then the transposed matrices
$(B_1^t,\ldots,B_D^t)$ define a covariant operator
with respect to $\tau$, too.
\par
If we omit the constant term in equation (\ref{defcovariant}),
we can be a little more general: In this situation we can also admit
matrices $B_j$ that are not quadratic, i.e.\ we can consider operators
${\cal D}_B : C^\infty(\R^D, \K^N) \to C^\infty(\R^D, \K^{M})$.
\begin{proposition}
Let $\tau$ be a representation of the group ${\bf G}$ in
${\rm Aut}\,\K^N$ and let $\sigma$ be a representation of ${\bf G}$ in
${\rm Aut}\,\K^M$.
Let $B_1,\ldots,B_D \in {\cal M}_{M\times N}$ and put
$B=(B_1,\ldots,B_D)$.
\par
We consider the operator ${\cal D}_B$ defined in equation (\ref{covoperator})
and put $m=0$.
Let $T_g$ denote the action of $\tau$ in
$C^{\infty}(\R^D,\K^N)$
and let $S_g$ denote the action of $\sigma$ in
$C^{\infty}(\R^D,\K^M)$.
The following statements are equivalent:
\begin{itemize}
\item[(a)]
The form of ${\cal D}_B$ does not change if we make a coordinate
transformation in
$\R^D : x\mapsto g\, x\ ,\ g\in{\bf G},$
and simultaneously coordinate transformations in
$\K^N : y \mapsto \tau (g)\, y $
and in
$\K^M : z \mapsto \sigma (g)\, z $.
\item[(b)]
${\cal D}_B$ intertwines $T_g$ and $S_g$ :
{ \renewcommand{\arraystretch}{2.5}
$$
\begin{array}{ccc}
  C^{\infty}(\R^D,\K^N) &\quad{\buildrel{\cal D}_B \over\longrightarrow}\quad &
            C^{\infty}(\R^D,\K^M)  \\
T_g \Big\downarrow & &\Big\downarrow  S_g \\
  C^{\infty}(\R^D,\K^N) & \quad {\buildrel {\cal D}_B \over \longrightarrow}
   \quad & C^{\infty}(\R^D,\K^M)
\end{array}
$$
}
\item[(c)]
$$
\sum_{k=1}^D g_{jk}\; \sigma (g)\; B_k\; \tau (g^{-1}) = B_j \qquad
\forall j\in\{1,\ldots,D\}\qquad \forall g \in {\bf G}
$$
where $g_{jk}$ are the components of $g\in {\bf G}$.
\end{itemize}
\end{proposition}
\par
For a given $\tau$ and $\sigma$ the set of all operators fulfilling
one of the conditions of Proposition $2.3$ will be denoted as
$Cov((\tau ,\K^{N});(\sigma,\K^{M}))$.
The following lemma is the infinitesimal version of the transformation
properties (\ref{defcovariant}).
\begin{lemma}
Let {\fraktur g} denote the Lie algebra of {\bf G}, and let
$L_{\alpha},\ \alpha\!\in\!\{ 1,\ldots,l\}$,
be a family of generators of
{\fraktur g}.
\par
A necessary condition that
a $D$-tuple of matrices
$B=(B_1,\ldots,B_D)$ defines a covariant operator
${\cal D}_B$
with respect to the representation $\tau$ is that
\begin{equation}
\label{infinitesimal}
\sum_{k=1}^D\;(L_{\alpha})_{jk}\;B_k =
[\, B_j\, , \, d\tau(L_{\alpha}) \, ]
\qquad\forall \alpha\in\{1,\ldots,l\}\quad
\forall j\in\{1,\ldots,D\}
\end{equation}
where $d\tau$ denotes the differential of $\tau$.
\par
If {\bf G} is connected, condition (\ref{infinitesimal}) is also
sufficient.
\end{lemma}
{\bf Sketch of the proof:}
\par
The infinitesimal form follows easily from the global condition.
Therefore we shall concentrate on the proof of the inverse implication.
First we show that the statement to be proven holds for one-parameter groups.
Let us take
the one-parameter group $g(t)=e^{itL_{\alpha}}$ and its representation
$T_{g}(t)=e^{itd_{\tau}(L_{\alpha})}$. By the commutator  expansion
$$\sum_{k} T_{g}B_{k}T_{g}^{-1}g_{ik}=\sum_{k}\sum_{n\geq 0}
\frac{i^{n}t^{n}}{n!}[d_{\tau}(L_{\alpha}),...,[d_{\tau}(L_{\alpha}),B_{k}]...]
(e^{itL_{\alpha}})_{ik}$$
and iterating
$$[d_{\tau}(L_{\alpha}),...,[d_{\tau}(L_{\alpha}),B_{k}]...]=
-[d_{\tau}(L_{\alpha}),...,[d_{\tau}(L_{\alpha}),
\sum_{k_{1}}B_{k_{1}}(L_{\alpha})_{kk_{1}})]...]$$
$$=-\sum_{k1}(L_{\alpha})_{kk_{1}}[d_{\tau}(L_{\alpha}),...,
[d_{\tau}(L_{\alpha}),B_{k1}]...]=(-1)^{n}\sum_{k_{1}...k_{n}}
(L_{\alpha})_{kk_{1}}(L_{\alpha})_{k_{1}k_{2}}...(L_{\alpha})_{k_{n-1}k_{n}}
B_{k_{n}}$$
so that
$$\sum_{k} T_{g}B_{k}T_{g}^{-1}g_{ki}=\sum_{n\geq 0}
\frac{(-i)^{n}t^{n}}{n!}\sum_{k}\sum_{k_{1}...k_{n}}(-1)^{n}(e^{itL_{\alpha}})
_{ik}
(L_{\alpha})_{kk_{1}}(L_{\alpha})_{k_{1}k_{2}}...(L_{\alpha})_{k_{n-1}k_{n}}
B_{k_{n}}$$
$$=\sum_{n\geq 0}
\frac{(-i)^{n}t^{n}}{n!}\sum_{k_{n}}(e^{itL_{\alpha}}L_{\alpha}^{n})_{ik_{n}}
B_{k_{n}}=
\sum_{k}[e^{itL_{\alpha}}(\sum_{n\geq 0}
\frac{(-i)^{n}t^{n}}{n!}L_{\alpha}^{n})]_{ik}B_{k_{n}}=B_{i}.$$
This ends the proof for one-parameter group.
Because
$$\sum_{k}T_{g_{2}g_{1}} B_{k} T_{g_{2}g_{1}}^{-1}(g_{2}g_{1})_{lk}=
\sum_{k}\sum_{i}T_{g_{2}}T_{g_{1}} B_{k} T_{g_{1}}^{-1}T_{g_{2}}^{-1}g_{1,ik}
g_{2,li}=
\sum_{i}T_{g_{2}}B_{i}T_{g_{2}}g_{2,li}=B_{l}$$
and the fact that the statement holds
for the one-parameter groups we have proven the
implication for group elements $g$ which are products
$g=g_{1}(t_{1})...g_{k}(t_{k})$ of elements from the one-parameter groups
$g_{i}(t_{i})$.
The set of such products is dense in some open subset $U$ containing
the identity.
By the continuity argument the global condition is
fulfilled on $U$, and consequently is fulfilled on the
connected component containing the identity.
$\Box$
\smallskip
\par
\begin{remark}
Let the Lie group ${\bf G}$ be the union of connected components
${\bf G}=\bigcup _{\alpha} {\bf G}^{\alpha}$ with ${\bf G}^{0}$ being
the connected component containing the unit element ${\bf e}$. Assume that
there exist(s) $R_{\alpha} \in {\bf G}$ such that
$R_{\alpha}{\bf G}^{0}={\bf G}^{\alpha}$. If for a given
representation $\tau$ equations (7) hold and if
\begin{equation}
\sum _{k=1}^{D} (R_{\alpha})_{jk}\tau (R_{\alpha})B_{k}\tau ^{-1}(R_{\alpha})=
B_{j}
\end{equation}
then the $D$-tuple $\{ B_{j} \} _{j=1,...,D}$ defines a
covariant operator ${\cal D}$ under the action of component(s)
${\cal G}^{\alpha}$.
\end{remark}
Similarly we can also prove:
\begin{lemma}
Let ${\bf G}$, $g$, $L_{\alpha}$ be as in Lemma $(2.5)$ and let $\sigma$,
$\tau$ be two representations of ${\bf G}$ in $\K^{N}$ and in
$\K^{M}$ respectively. A necessary condition that a $D$-tuple of matrices
$B=(B_{1},...,B_{D}
)$ defines a covariant operator ${\cal D}_{B}
\in Cov((\tau ,\K^{M});(\sigma ,\K^{N}))$ is that:
\begin{equation}
\sum _{k=1}^{D}(L_{\alpha})_{jk}B_{k}+d\sigma (L_{\alpha})B_{j}+B_{j}d\tau
(L_{\alpha})=0
\end{equation}
for all $j,k \in \{ 1,...,D \}$ and $\alpha =1,...,dim{\bf G}$.\par
If ${\bf G}$ is connected this condition is also sufficient.
\end{lemma}
\smallskip
For the case of the rotation group
${\bf SO(3)}$ in three-dimensional space and also for the
proper orthochronous Lorentz group
${\bf L}^{\uparrow}_{+}{\bf (4)}$
in four-dimensional space-time
covariant operators have been extensively studied,
cf.\ \cite{GelfandShapiro,Ljubarski,WaveEquations}
and the references therein.
\par
In the sequel we want to study the inverse of a given covariant operator.
It is therefore natural to ask whether we can find any elliptic operators
in ${\rm Cov}(\K^N,\tau)$.
\par
For an operator
${\cal D}_B = \sum_{j=1}^D B_j {\partial \over \partial x_j} + m E$
and a differential form
$\sum_{j=1}^D p_j\, dx_j$
we define the characteristic
polynomial in the usual way:
$$\sigma_{{\cal D}_B}(p_1,\ldots,p_D)\;
{\buildrel {\rm def}\over =}\;
i\;\sum_{j=1}^D B_j\, p_j \ \ .
$$
Note that this definition depends in general on the choice of coordinates.
\begin{lemma}
\label{charpolynomial}
${ }$
\par
\begin{itemize}
\item[(a)]
Let ${\bf G}\subseteq{\bf O(D)}$ and let
${\cal D}_B \in {\rm Cov}(\K^N,\tau)$.
\par
The form of $\sigma_{{\cal D}_B} (p_1,\ldots,p_D)$
does not change if we make a coordinate transformation
in $\R^D$ :
$x\mapsto g x $, $g\in {\bf G}$,
and simultaneously a coordinate transformation in $\K^N$:
$y \mapsto \tau(g)\, y$.
\item[(b)]
Let ${\bf G}$ be either ${\bf SO(D)}$ or ${\bf O(D)}$ and let
${\cal D}_B \in {\rm Cov}(\K^N,\tau)$.
\par
We have
\begin{equation}
{\rm det}\bigl(\sigma_{{\cal D}_B}(p_1,\ldots,p_D)\bigr) =
C\,(p_1^2 + \ldots + p_D^2)^{n}
\end{equation}
for some constant $C \in \C$ and $n \in {\bf N}$.
\par
Moreover, if $N$ is odd, then
$C=0$, i.e.\  elliptic operators that are covariant with respect to some
representation of ${\bf SO(D)}$ or ${\bf O(D)}$
can only exist if the dimension of the
representation space is even.
\end{itemize}
\end{lemma}
\par
{\bf Proof:}
(a) is easily seen by employing the covariance condition
(\ref{defcovariant}).
To prove (b), observe that
${\rm det}\bigl(\sigma_{{\cal D}_B} (p_1,\ldots,p_D)\bigr)$
is invariant under rotations and must therefore be a function of
$p_1^2 + \ldots + p_D^2$.
The assertion now follows from the fact that
${\rm det}\bigl(\sigma_{{\cal D}_B} (p_1,\ldots,p_D)\bigr)$ must be
a polynomial and a homogeneous function of order less or equal to $N$.
\hfill $\Box$
\par\medskip
\begin{remark}
Let ${\bf G}$ be either ${\bf SO(D)}$ or ${\bf O(D)}$ and let $N$ be even.
For a covariant operator ${\cal D}_B \in {\rm Cov}(\K^N,\tau)$ we have
\begin{equation}
{\rm det} \Bigl( i \sum_{j=1}^D B_j\, p_j + m E \Bigr) =
C \prod_{\alpha=1}^{{N\over 2}} \bigl( p_1^2 + \ldots + p_D^2
 + r_{\alpha}^2 \bigr)
\end{equation}
where $r_{\alpha},\ \alpha = 1,\ldots,{N\over 2},$
and $C \neq 0$ are constants $\in \C$.
\par
If all $r_{\alpha}$ are real, the operator ${\cal D}_B$ is invertible
on suitably chosen function spaces and in this case we shall call it
admissible.
\end{remark}
\par
\smallskip
Given two different but equivalent representations $\tau$
and $\tilde\tau$, the following remark shows how we can identify
${\rm Cov}(\K^N,\tau)$ and
${\rm Cov}(\K^N,\tilde\tau)$.
\begin{remark}
\label{equivrepr}
We assume that
$B=(B_1,\ldots,B_D)$ defines a covariant operator
with respect to the representation $\tau$.
Let $\tilde\tau$ be an equivalent
representation:
$\tilde \tau(g) = M\,\tau(g)\, M^{-1}$.
\par
Then
$B^\prime = (B_1^\prime,\ldots,B_D^\prime),\ \
B_j^\prime =M\, B_j\, M^{-1}$,
defines a covariant operator with respect to $\tilde \tau$.
\end{remark}
\smallskip
\begin{remark}
It is possible to covariant differential operators of higher order, too.
For this let:
$$
{\cal D}_{n}=\sum _{\alpha : |\alpha| \leq n}  {\cal B}_{\alpha}\partial _
{\alpha } + {\cal M}
$$
where $\alpha = (\alpha _{1},...,\alpha _{D})$, $\alpha _{i} \in \N \cup \{
0 \}$, $|\alpha | = \alpha _{1}+...+\alpha _{D}$, ${\cal B}_{\alpha} \in
{\cal M}_{N \times N}(\K)$, $\partial _{\alpha}=\frac{\partial ^{\alpha
_{1}+...
\alpha _{D}}}{\partial x_{D}^{\alpha _{D}}...\partial x_{1}^{\alpha _{1}}}$ and
let $\tau$ be a representation of a group ${\bf G}$ in $\K^{N}$. Then the
operator ${\cal D}_{n}$ is called $\tau$-covariant differential operator of
order $n$ iff
\begin{description}
\item{(i)}  there exists ${\alpha}$ such that $|\alpha |=n$ and ${\cal
B}_{\alpha}
\neq 0$
\item{(ii)}
the following diagrams commute:
$$
\begin{array}{ccc}
  C^{\infty}(\R^D,\K^N) &\quad{\buildrel{\cal D}_{n} \over\longrightarrow}\quad
&
            C^{\infty}(\R^D,\K^N)  \\
T_{g}^{\tau} \Big\downarrow & &\Big\downarrow  T_{g}^{\tau} \\
  C^{\infty}(\R^D,\K^N) & \quad {\buildrel {\cal D}_B \over \longrightarrow}
   \quad & C^{\infty}(\R^D,\K^N)
\end{array}
$$
\end{description}
In particular, taking ${\cal D}^{1},...,{\cal D}^{n} \in Cov(\tau ;\K^{N})$
the operator ${\cal D}_{n}={\cal D}^{n}...{\cal D}^{1}$ is a covariant
operator of $n$-th order. However, since by increasing the dimension $N$ of
the target space $\K^{N}$ the $n$-order covariant equation can be reduced to
first order we shall mainly restrict ourselves to the first order operators.
\end{remark}
\smallskip
\par
Let us now focus on the case ${\bf G}={\bf SO(D)}$.
The representation theory of {\bf SO(D)} is well known, cf.\
\cite{Boerner,Broecker,Fulton,GelfandShapiro}.
An important question for physics is which representations $\tau$ of
${\bf G}={\bf SO(D)}$ admit an extension
to a representation $\tilde\tau$ of {\bf O(D)}.
Since {\bf SO(D)} is a subgroup of index 2 of {\bf O(D)}, it is a
normal subgroup and
${\bf O(D)}/{\bf SO(D)} \cong \Z_2$.
Taking any $M\in {\bf O(D)}\setminus {\bf SO(D)}$ it is easy to check
that $\tau$ can be extended to {\bf O(D)} iff there exists
$\tilde\tau(M) \in {\cal M}_{N\times N}(\K)$ such that
\begin{equation}
\label{extensioncondition}
\tau(M \cdot A \cdot M) = \tilde\tau(M)\cdot \tau(A) \cdot \tilde\tau(M)
\qquad \forall A \in {\bf SO(D)}\ \ .
\end{equation}
If $D$ is odd one can always extend a given representation $\tau$:
The fact that $D$ is odd implies that
the matrix
$ M = -E_D =(-\delta_{ij})$
has determinant $-1$, and if we put
$\tilde\tau(M) = \pm{\rm id}_V$, condition
(\ref{extensioncondition}) is fulfilled.
\par
Let us now have a look at
\begin{equation}
\label{defR}
 R=\left(\begin{array}{rl}
-1 & 0 \\
0  & E_{D-1} \\
\end{array}\right)
\end{equation}
which is the reflection at the hyperplane $\{x_1=0\}$.
The choice
$\tilde \tau(-E_D) = \pm {\rm id}_V$
implies that the reflection $R$ is represented by
\begin{equation}
\label{tildetauR}
\tilde\tau(R) =
\pm \tau \left(
\begin{array}{cc}
 1 & 0 \\
 0 & -E_{D-1} \\
\end{array} \right)\ \ .
\end{equation}
\par
The case of even dimension
is more complicated so that we only give a summary of some
group-theoretic results, referring the reader to \cite{Boerner} for
details.
\par
We assume that $\tau$ is an irreducible unitary representation. Taking some
$M\in{\bf O(D)}\setminus{\bf SO(D)}$, we consider the representation
$\sigma(A) = \tau(\,M^{-1}A\,M) , \  A\in{\bf SO(D)}$.
If $\sigma$ and $\tau$ are equivalent, $\tau$ is called
self-conjugate.
In this case $\tau$ can be extended to {\bf O(D)}, and the extension
is unique up to sign.
If, however, $\sigma$ and $\tau$ are not equivalent, one has to pass
to the induced representation $\tau_{\rm ind}$ of {\bf O(D)},
i.e. one has to
double the dimension of the representation space $\K^N$.
$\tau_{\rm ind}$ is an irreducible representation of {\bf O(D)}, and it
is the only irreducible representation of {\bf O(D)} which contains
$\tau$ when being restricted to {\bf SO(D)}.
\pagebreak[2]
\medskip
\par
Now we can introduce reflections into the concept of covariant
operators.
\begin{definition}
Let $\tau$ be a representation of ${\bf G}={\bf SO(D)}$,
and let $\tilde\tau$ be an extension
of $\tau$ to {\bf O(D)}.
\par
We call an operator ${\cal D}_B \in
{\rm Cov}(\K^N,\tau)$ reflection-covariant with respect to $\tilde\tau$ iff
it transforms covariantly under the full orthogonal group, i.e. if
(\ref{defcovariant})
holds $\forall g\in{\bf O(D)}$.
\end{definition}
\begin{remark}
Let ${\cal D}_B$ be a covariant operator with respect to a representation
$\tau$ of {\bf SO(D)}, and let $\tilde\tau$ be an extension of $\tau$
to {\bf O(D)}.
\smallskip
${\cal D}_B$ is reflection-covariant with respect to $\tilde\tau$
if and only if
\begin{eqnarray}
\tilde\tau(R)\,B_1\,\tilde\tau(R) &=&- B_1\nonumber  \\
\tilde\tau(R)\,B_j\,\tilde\tau(R) &=&B_j \qquad \forall j \in\{2,\ldots, D\}
\end{eqnarray}
\par
where $R$ is the matrix in equation (\ref{defR}).
\end{remark}
\medskip
\par
Unitary representations of the classical groups are well understood.
In the sequel we have to  use representations in terms of real matrices.
\par
Let $V$ be a complex finite-dimensional vector space.
Given a representation
$\tau\ :\ {\bf G} \to {\rm Aut}V$,
it is natural to ask whether
$\tau$ can somehow be transformed into a representation in terms of
real matrices.
A comprehensive treatment of this question can be found in
\cite{Broecker,Fulton}.
\par
\pagebreak[2]
$\tau$ is {\it of real type} iff there is an antilinear map
$J\ :\ V\to V$ such that
$J^2={\rm id}_V$ and
$J \tau(g) = \tau(g) J\quad \forall g\in{\bf G}$.
\par
If $\tau$ is of real type, consider
$W=\{ x\in V\ \vert\ x=Jx \}$.
$W$ is a real subspace which is $\tau (g)$-invariant
$\forall g\in{\bf G}$.
We have the decomposition $V=W\oplus i W$ which shows that $\tau$ can
be obtained from $\tau_{\rm real} : {\bf G} \to W$ by extending the field
of scalars.
Choosing a basis for $W$, we get a representation in terms of real
matrices.
\par
$\tau$ is {\it of quaternionic type} iff there is an antilinear map
$J : V\to V$ such that $J^2 = -{\rm id}_V$ and
$J \tau(g) = \tau(g) J\quad \forall g\in{\bf G}$.
If the representation $\tau$ is of quaternionic type, it can be extended
to
$\tau_{\rm quat} : V\oplus j V$, where $\{ 1,i,j,k\}$ denotes, as usual,
the canonical basis for the space of quaternions.
\par
If $\tau$ is neither of real nor of quaternionic type, we say that $\tau$
is {\it of complex type}.
The following proposition is a well-known criterion
to determine the type of a given irreducible representation.
\begin{proposition}
Let $dg$ denote the normalized Haar measure on the compact Lie group
{\bf G} and let
$\chi_{\tau}$ denote the character of the irreducible representation
$\tau : {\bf G} \to {\rm End}V$.
$$
\int_{\bf G} \chi_{\tau}(g^2)\; dg =
\left\{ \begin{array}{rcl}
1 & \iff &
\tau\ {\it is\ of\ real\ type} \\
0 & \iff &
\tau\ {\it is\ of\ complex\ type} \\
-1 & \iff &
\tau\ {\it is\ of\ quaternionic\ type}\\
\end{array}
\right.
$$
\end{proposition}
\subsection{Non-Gaussian Noise}
In this section we shall deal with {\bf G}-invariant and
reflection-positive noise.
Since mathematical physicists might be less acquainted with
the notion of non-Gaussian noise, we briefly review some basic
definitions and facts.
\begin{definition}
Let $(\Omega, \Sigma, \mu)$ be a probability space, and
let $T$ be a space of smooth test functions
$\R^D \to \R^N$.
We assume that $T$ is equipped with some topology.
\par
A generalized random field indexed by $T$ is a map
$$
\varphi \ : \ T\ \ \longrightarrow\ \
\big\{ {\rm real\!-\!valued\ random\ variables\ on\ } \Omega \big\}
$$
which is linear almost surely, i.e.
$\forall f,g \in T,\ \forall \lambda \in \R$
\begin{eqnarray*}
   \varphi (f+g)&=&\varphi (f) + \varphi (g) \\
   \varphi (\lambda f) & =& \lambda \varphi(f) \ \
\end{eqnarray*}
and which is continuous in the sense that if
$f_n \to f$ in $T$ then $\varphi(f_n) \to \varphi(f)$ in
probability.
\end{definition}
On the formal level, we have
$$ \varphi(f) = \langle \varphi , f \rangle
       =\sum_{\alpha =1}^N\;\langle \varphi_{\alpha} , f_{\alpha} \rangle
       =\sum_{\alpha =1}^N\ \int_{\R^D}\; \varphi_{\alpha}(x) f_{\alpha}(x)\,
d\!x
       \ \ .
$$
\pagebreak[2]
\begin{definition}
\hfill\break
Let
${\cal D}\!=\!{\cal D}(\R^D)\otimes \R^N$
denote the space of
$C^{\infty}$-functions $\R^D \to \R^N$ with compact support.
White noise is a generalized random field $\varphi$ indexed
by ${\cal D}$
such that its characteristic functional is given by
\begin{equation}
\label{noisefunctional}
\Gamma(f) =
E\Bigl( e^{i \varphi (f) } \Bigr)
 =
e^{- \int_{\R^D} \psi(f(x))\;dx}\ \ \ .
\end{equation}
The function $\psi : \R^N \to \C$ has the so-called L\'evy-Khinchin
representation
{\hfuzz=20pt
\begin{equation}
\label{LevyKhinchin}
\psi(y) = i <\beta, y> + \frac{1}{2}< y, Ay> +\!\!\!\int
\limits_{\R^N \setminus \{ 0\} }
\!\!\!\Bigl( 1 - e^{i <\alpha , y> } + {{ i < \alpha, y>} \over
{1 + \Vert \alpha \Vert^2} } \Bigr)\;{ {1 + \Vert \alpha \Vert^2} \over
\Vert \alpha \Vert^2  }\; d\!\kappa(\alpha)
\end{equation}}
\par
where $\beta\in\R^N$, $A$ is a non-negative definite $N\times N$-matrix
and $\kappa$ is a non-negative, bounded measure on $\R^N \setminus \{ 0\}$.
\par
If $\kappa=0$ and $A\ne 0$, $\varphi$ is called Gaussian white noise whereas
in the case $A=0$, $\kappa\ne 0$, $\varphi$ is called Poisson noise.
In the following we put always $\beta =0$ for simplicity.
\end{definition}
In the last section we mentioned that we need
representations of the group {\bf G} in terms of
real matrices.
The reason for this is that
${\cal D}(\R^D)\otimes\R^N$
is a vector space over $\R$.
\par
Since ${\cal D}$ is a nuclear space, by Minlos\rq\ theorem
(cf.\cite{GelfandVilenkin})
there is a unique
probability measure $\mu$ on the dual space ${\cal D}^{\prime}$
such that
$$
\int_{{\cal D}^{\prime}} e^{i (\eta , f)} d\!\mu (\eta) = \Gamma(f)
$$
where $(\,\cdot\, , \, \cdot \,)$ denotes the canonical
pairing between ${\cal D}^{\prime}$ and ${\cal D}$.
\par
The function $\psi$ in (\ref{LevyKhinchin}) is a negative definite
function, cf.\ \cite{BergForst}.
$$
\psi_G(y) = \frac{1}{2} <y , Ay>
$$
is the Gaussian part and
$$
\psi_P(y) = \int_{\R^N \setminus \{ 0\} }
\Bigl( 1 - e^{i <\alpha , y> } + {{ i < \alpha, y>} \over
{1 + \Vert \alpha \Vert^2} } \Bigr)\;{ {1 + \Vert \alpha \Vert^2} \over
\Vert \alpha \Vert^2  }\; d\!\kappa(\alpha)
$$
is the Poisson part of $\psi$.
\par
We shall also use the notation
$$
\Gamma_G (f) = E_G (e^{i\, \varphi_G(f)})
= e^{-\int \psi_G (f(x)) dx }
$$
and the analogous notation for the Poisson part.
\par
The noise $\varphi$ can be regarded as the sum of Gaussian and Poisson
noise:
$\varphi = \varphi_G + \varphi_P$.
Correspondingly, we have a measure $\mu_G$ and a measure $\mu_P$ on
${\cal D}^{\prime}$, and $\mu$ is the convolution of these two measures:
$\mu = \mu_G \ast \mu_P$.
\par
Let us mention two characteristic properties of white noise.
White noise is invariant under translations in the sense that
the random variables $\varphi(f_{x_0})$ and $\varphi(f)$ are
equal in law, where $f_{x_0}$ is the function
$x\mapsto f(x+x_0)$.
\par
If we take two functions $f_1, f_2 \in{\cal D}$ with disjoint
supports, the random variables $\varphi(f_1)$ and $\varphi(f_2)$
are independent.
\smallskip
\par
If $\varphi$ is white noise such that the random variables
$\varphi(f)$ have zero mean and finite second moments $\forall f$,
the function $\psi$ in (\ref{LevyKhinchin}) has the so-called
Kolmogorov canonical representation
\begin{equation}
\label{Kolmogorovrepresentation}
\psi(y) =\,\frac{1}{2}<y,Ay>\; +
\int\limits_{\R^N \setminus\{ 0\} } \Bigl( 1 -
e^{i <\alpha, y>} + i <\alpha , y> \Bigr) d\!\nu(\alpha)
\end{equation}
where the so-called L\'evy measure $\nu$ has the property
$\int_{\R^N \setminus \{ 0\} } \Vert \alpha\Vert^2 d\!\nu(\alpha) < \infty$.
In this case $\psi$ satisfies the inequality
$\vert \psi(y)\vert \le M \Vert y\Vert^2\quad \forall y\in \R^N$
where $M$ is some constant $\ge 0$.
This makes it possible to extend the generalized random field $\varphi$
to $L^2$.
\par
In the sequel we shall restrict the class of admitted characteristic
functionals even further.
We shall assume that the L\'evy measure $\nu$ in
(\ref{Kolmogorovrepresentation}) is invariant under the reflection
$\alpha\mapsto - \alpha$.
Under this assumption the characteristic functional corresponding to the
Poisson part is of the form
\begin{equation}
\label{ournoise}
\Gamma _{P}(f) = E_{P}(e^{i \varphi(f) } )
=e^{ \int_{\R^D} \int_{\R^N} \bigl( e^{ i < \alpha, f(x)>} - 1 \bigr)
d\!\nu(\alpha)\;dx }\ .
\end{equation}
Moreover, we assume that the measure $\nu$ satisfies the condition
\begin{equation}
\int_{\R^N} e^{ t\;\Vert \alpha\Vert } d\!\nu(\alpha)
<\infty  \qquad
\forall t \ge 0\ \ \ .
\end{equation}
\par
This condition guarantees the existence of all moments of the corresponding
noise and, moreover, it allows us to extend the characteristic
functional as an analytic function.
To be more precise, for fixed $f\in {\cal D}(\R^{D}) \otimes \R^{N}$
\begin{equation}
\C^{N} \ni \xi \longmapsto \Gamma _{P}(\xi f) =
E _{P}\bigl(e^{i\,(\,\cdot\, ,\xi f)}\bigr)
\end{equation}
is an entire function in $\xi$ obeying the estimate:
\begin{equation}
\vert\Gamma _P (\xi F)\vert \le \exp \Bigl( \vert\xi\vert \int \Vert f(x)
\Vert \;d\!x \cdot
\int \Vert\alpha\Vert\; e^{(\,\vert\xi\vert\;\,\Vert\alpha\Vert\;\,
\Vert |f|\Vert\,)}
d\!\nu(\alpha) \Bigr)
\end{equation}
where
$\Vert f(x)\Vert =\Bigl( \,\sum\limits_{i=1}^N \vert f_i (x)\vert^2
\Bigr)^{\frac{1}{2}}$ and $\Vert |f|\Vert = \sup_{x}|f(x)|$.
\begin{lemma}
Let us assume that the L\'evy measure $\nu$ in (\ref{ournoise}) has finite
first
order moments.
Then for any $f\in D(\R^{D})$, any cylinder
function $F\in L^2(\mu_P)$ which is bounded and $C^1$
the following integration by parts formula
holds:
{\hfuzz=5pt
$$
\int\limits_{{\cal D}^{\prime}(\R^{D}) \otimes \R^{N}}\!\!\!<\eta ,f^{\lambda
}>
F(\eta)\;
d\!\mu(\eta) =
$$
\begin{equation}
\label{intbyparts}
\int\! <f^{\lambda}(x),E\bigl( A \frac{\delta}{\delta \eta(x)}
F(\eta)\bigr)>d\!x
\;+\;
\int
\int f^{\lambda}(x)\; E\bigl( F(\eta +\alpha \delta
(x- \cdot )\bigr)\;\alpha_\lambda\;d\!\nu(\alpha) \;d\!x
\end{equation}
}
where
$(f^{\lambda})_{i}=\delta ^{\lambda}_{i}f$\ ,
$\frac{\delta}{\delta \eta (x)}$ denotes the functional derivative
(widely used in mathematical physics see e.g. \cite{Glimm-Jaffe}, and
$\bigl(A\frac{\delta}{\delta
\eta (x)}\bigr)_j=\sum\limits_{k=1}^{N} A_{jk}\frac{\delta}
{\delta \eta_k (x)}$\ .
\end{lemma}
{\bf Proof:}\\
Take $F(\eta)=\exp{i(\eta, g)}$.
Employing (\ref{ournoise}), it is easily seen that (\ref{intbyparts})
holds.
Since any bounded $C^{1}$ cylinder
function can be uniformly approximated by the sums $\sum_{n}c_{n}\exp{i(\eta  ,
g)}$ (see e.g. \cite{Glimm-Jaffe}), the assertion follows.\hfill $\Box$
\par
\bigskip
\par
If the characteristic functional of Poisson noise $\varphi$ is of the form
(\ref{ournoise}), the moments of $\varphi$ are given by
\begin{equation}
\label{Poissonmoments}
E_{P}\bigl(\,\prod _{i=1}^{n}(\varphi ,f_{i})\bigr)
= \sum_{\parbox[t]{60pt}{
$\scriptstyle \Pi _{1}\cup ...
\cup \Pi _{k}=J_{n}$ \\
$\scriptstyle \Pi _{\alpha}\cap \Pi _{\beta}=\emptyset$ \\
$\scriptstyle {\rm for}\; \alpha \neq \beta$} }
\prod _{l=1}^{k} \int \int
\prod _{j \in \Pi _{l}}<\alpha_l ,
f_j (x_l)>\;d\!x_l\;d\!\nu(\alpha_l)
\end{equation}
where the summation runs over the set of all partitions of
$J_{n}=\{ 1,2,...,n \}$.
If the noise $\varphi$ is the sum of a Gaussian and a Poisson part,
formula (\ref{Poissonmoments}) has to be altered:
\begin{equation}
\label{moments}
E\bigl(\,\prod _{i=1}^{n}(\varphi ,f_{i})\bigr)=\sum_{\parbox[t]{50pt}{
$\scriptstyle \Pi_G \cup \Pi_P=J_{n}$ \\
$\scriptstyle \Pi_G \cap \Pi_P = \emptyset $} }
E_P\bigl(\prod _{i \in \Pi_P}(\varphi ,f_{i})\bigr)
E_G\bigl(\prod _{i \in \Pi_G}(\varphi ,f_i)\bigr)
\end{equation}
The moments of the Gaussian part are uniquely determined by the
covariance $A$:
\begin{equation}
E_G\bigl((\varphi ,f_1)(\varphi ,f_2)\bigr)=\int<f_1(x),A\,f_2(x) > d\!x.
\end{equation}
\begin{remark}
The number of terms in (\ref{Poissonmoments}) is
$\Pi _n=\sum_{p=1}^n S_n (p)$, where $S_n (p)$ are the so-called
Stirling numbers of the second kind.
They are given explicitly by
\begin{equation}
S_n (p)=\frac{1}{p!}
\sum_{j=0}^p (-1)^j {p \choose j}(p-j)^k\ .
\end{equation}
Therefore the total number of terms in (\ref{moments})
is $\sum_{k=0}^n{n \choose  k} \Pi _k$.
\end{remark}
To derive (\ref{ournoise}) we made the assumption that the L\'evy measure $\nu$
is
invariant under the reflection
$\alpha \mapsto - \alpha$.
This implies that the contributions coming from
partitions $(\Pi _{\alpha})$ in (\ref{Poissonmoments}) containing some
$\Pi_\alpha$ with an odd number of elements vanish.
\par
\smallskip
The following remark shows that the carrier set of Poisson noise is
extremely small: it consists of locally finite linear combinations of
delta distributions.
\begin{remark}
\label{carrierProp}
Let
$$
C_{lf}(\R^D)= \{\Lambda\subset\R^D\ \vert\ \Lambda\cap K {\rm\ is\ finite\ for\
every\ compact\ set\ } K\}\ ,
$$
i.e.\ $C_{lf}$ is the set of \lq locally finite configurations\rq .
$C_{lf}$ can be given a topology such that $C_{lf}$ is a complete metrizable
space, cf.\ \cite{Kerstan}.
\par
If $\Lambda\in C_{lf}(\R^D)$, $\Lambda$  obviously contains either a
finite number of points or countably many points. Let us fix an
enumeration of these points,
i.e.\
$\Lambda= (x_1, x_2,\ldots)$, $\ x_i \in \R^D$.
Take $\Gamma=(\gamma_1, \gamma_2,\ldots )\in (\R^{N})^{\vert\Lambda\vert}$.
We define
\begin{equation}
\delta(\Lambda,\Gamma)(x) = \sum_{{x_i\in\Lambda}\atop{\gamma_i\in \Gamma}}
\gamma_i\; \delta(x - x_i)
\end{equation}
where $\forall f=(f_1,\ldots,f_N)\in {\cal D}(\R^D)\otimes\R^N$
\begin{equation}
\bigl( \gamma_i\; \delta(\,\cdot\, - x_i)\, ,\, f\bigr)
= \sum_{k=1}^N (\gamma_i)_k\; f_k(x_i)\ .
\end{equation}
Adapting the argument in \cite{Kerstan}, it can be proved that the set
$$
{\cal C} = \bigl\{ \delta(\Lambda,\Gamma)\ \vert\ \Lambda \in C_{lf}(\R^D) ,
 \Gamma\in ({\rm supp}\; \nu)^{\vert \Lambda\vert} \bigr\}
$$
is a carrier set for $\mu_P$,
i.e.\ $\mu_P({\cal C})=1$.
\end{remark}
\begin{remark}
Let $\Lambda$  be an open subset of $\R^{D}$. We define $\sigma $-algebra
$\Sigma (\Lambda)$ as a minimal ($\mu_{P}$-complete) $\sigma $-algebra
of sets generated by random elements $(\varphi ,f)$ with $f \in
{\cal D}(\R^{D})$
supported in $\Lambda$. For $\Lambda$ closed we define $\Sigma (\Lambda)$ as an
intersection of all $\Sigma (\Lambda ')$, where $\Lambda '$ is open and
$\Lambda
\subset \Lambda '$. Let $\Gamma \subset \R^{D}$ be a closed subset of $\R^{D}$
and of (Lebesgue)measure zero. It can be easily deduced from Remark
\ref{carrierProp}
that then $\Sigma (\Gamma)$ is a trivial $\sigma $-algebra. From this it
follows that the random field $\mu_{P}$ has a Markov property in the
following sense:
\par
for any open $\Lambda \subset \R^{D}$ with sufficiently regular boundary
$\partial \Lambda$ and any bounded $F$,$G$ measurable with respect
$\Sigma (\Lambda)$ respectively $\Sigma (\Lambda ^{c})$:
\begin{equation}
E_{\mu_{P}}\{ F\cdot G|\Sigma(\partial \Lambda)\} =
E_{\mu_{P}}\{ F|\Sigma(\partial \Lambda)\} \cdot
E_{\mu_{P}}\{ G|\Sigma(\partial \Lambda)\}=
E_{\mu_{P}}(F)\cdot E_{\mu_{P}}(G),
\end{equation}
where $E_{\mu_{P}}\{ -|\Sigma(\cdot)\} $ denotes the corresponding conditional
expectation value of $(-)$ with respect to the $\sigma $-algebra
$\Sigma(\cdot)$.
\end{remark}
\begin{definition}
Let $\tau$ be a representation of the group $(S)O(D)$ in the space $\R^{N}$.
We will say that the random field $\varphi$ (given by (\ref{noisefunctional})
and (\ref{LevyKhinchin})) is
$\tau$-covariant random field iff
\begin{equation}
E(e^{i(\varphi,T_{\tau}f)})=E(e^{i(\varphi,f)})=
E(e^{i(T_{\tau}^{\ast}\varphi,f)})
\end{equation}
for all $f\in {\cal D}(\R^{D})\otimes \R^{N}$. $T^{\ast}_{\tau}$
means the adjoint of the representation $T_{\tau}$ acting in the space
${\cal D}'(\R^{D})\otimes \R^{N}$ under the canonical pairing
${}_{{\cal D}'}(~,~)_{{\cal D}}$.
\end{definition}
\begin{lemma}

Let $\tau$ be a representation of $(S)O(D)$ in the space $\R^{N}$ and let
$\varphi$ be a white noise given by (\ref{noisefunctional}) and
(\ref{LevyKhinchin}). Then the noise $\varphi$ is
$\tau$-covariant iff
\par
(i)  $\beta =0$
\par
(ii) $\tau ^{T}A\tau =A$
\par
and
\par
(iii)the measure $d\kappa$ is $\tau $-invariant
(providing that $\tau$ is given by orthogonal matrices).
\end{lemma}

Let $\varphi$ be $\tau$-covariant white noise and let $R_{\tau}$ be the
representative of the reflection operator $R$ for the representation
$\tau$ (see \ref{extensioncondition}). Let $f^{\alpha}\in {\cal D}(\R^{D})
\otimes \R^{N}$ be
a finite sequence of test functions supported on $\{ (t,x) \in \R^{D}|t>0 \}$.
Then for any finite sequence $c_{\alpha} \in \C $ we have:
\begin{equation}
\sum _{\alpha ,\beta}c_{\alpha}\overline{c_{\beta}}\Gamma (
e^{i(\varphi ,f^{\alpha})}e^{-i(\varphi ,R_{\tau}f^{\beta})})
=\sum _{\alpha ,\beta}c_{\alpha}\overline{c_{\beta}}\Gamma (
e^{i(\varphi ,f^{\alpha})})\Gamma (e^{-i(\varphi ,R_{\tau}f^{\beta})})
=|\sum _{\alpha }c_{\alpha}\Gamma (e^{i(\varphi ,f^{\alpha})})|^{2} \geq 0
\end{equation}
providing that the noise is $R_{\tau}$-invariant.
\begin{remark}

The last property express the so called reflection positivity of the noise
$\varphi $. Taking such reflection positive and covariant noise one can
construct from the
moments of it (see e.g. \cite{Glimm-Jaffe},\cite{Simon}) some covariant
quantum field fulfilling all
Wightman axioms. However, it is fairly easy to show that the arising quantum
field theory is a multiple of the identity operator.
\end{remark}

\subsection{Covariant SPDEs and Their Solutions}

Let ${\cal D} \in Cov(\tau ,\R^{N})$ for some real representation
$\tau$ of $SO(D)$ and let $\tilde{{\cal D}}$ be the transpose of
${\cal D}$ in the canonical pairing ${}_{{\cal S}'(\R^{D})\otimes \R^{N}}
<\cdot ~,~\cdot >_{{\cal S}(\R^{D})\otimes \R^{N}}$.
We shall consider stochastic partial differential equation (SPDE) of the
type
\begin{equation}
\label{spde}
\tilde{{\cal D}}\varphi =\eta
\end{equation}
where $\eta$ is given generalized random field indexed by
${\cal S}(\R^{D})\otimes \R^{N}$. An operator ${\cal D}$ will
be called regular (corresp. the equation will be called regular) iff
there exists a nuclear space ${\cal F}$ such that the (principial)
Green function ${\cal D}^{-1}$ of ${\cal D}$ is defined on ${\cal F}$
and ${\cal D}^{-1}$ maps ${\cal F}$ continuously into
${\cal S}(\R^{D})\otimes \R^{N}$. A generalized random field
$\varphi$ indexed by ${\cal F}$ is called a weak solution of regular
equation (\ref{spde}) iff:
\begin{equation}
\label{spdesolu}
<\varphi , f> \cong <\eta ,{\cal D}^{-1}f>~~~\hbox{for all}~f \in {\cal F}
\end{equation}
where $\cong$ means equality in law. Denoting by
$\Gamma _{\eta}$ the characteristic functional of the field $\eta$ we
have that the characteristic functional $\Gamma _{\varphi}$ of a weak
solution $\varphi$ of regular equation (\ref{spde}) is given by:
\begin{equation}
\Gamma_{\varphi}(f)=\Gamma_{\eta}({\cal D}^{-1}f)~~~\hbox{for}~f \in
{\cal F}
\end{equation}

The case in which ${\cal D}:{\cal S}(\R^{D})\otimes \R^{N} \rightarrow
{\cal S}(\R^{D})\otimes \R^{N}$ is continuous bijection will be called
strongly regular. For example, if ${\cal D}$ is admissible with strictly
positive mass spectrum then ${\cal D}$ is strongly regular.
In the case of strongly regular ${\cal D}$ the space ${\cal F}(\R^{D})$
could be chosen as ${\cal S}(\R^{D})\otimes \R^{N}$.

Let ${\cal K}_{{\cal D}} \equiv \{ \chi \in {\cal S}'(\R^{D})\otimes \R^{N}
|~\tilde{{\cal D}}\chi =0 \}$. Then for any weak solution $\varphi $ of a
regular equation (\ref{spde}) and for any $\chi \in {\cal K}_{{\cal D}} \cap
{\cal F}'$ the new random field $\varphi _{\chi}$ the characteristic functional
of which is given by:
\begin{equation}
\Gamma _{\varphi _{\chi}}(f)=e^{i<\chi ,f>}\Gamma _{\varphi}(f)
\end{equation}
is again a weak solution of (\ref{spde}). In fact it could be proven, that
fixing
the space ${\cal F}$ the whole set of weak solutions of (\ref{spde}) could be
exhausted in this way.

Let us recall that a generalized random field $\eta$ indexed by a space
${\cal F}$ is called $\tau$-covariant iff
\par
(i)  $T^{\tau}_{g}$ acts in the space ${\cal F}$,
(ii) $<\eta ,T^{\tau}_{g}f> \cong <\eta ,f>$ for each $g \in (S)O(D)$ and
$f \in {\cal F}$.

\begin{proposition}
Let us consider a regular equation (\ref{spde}) with $\eta$ being
$\tau$-covariant.
Then the weak solution of (\ref{spde}) given by (\ref{spdesolu}) is again
$\tau$-covariant
random field (providing $T^{\tau}_{g}$ acts in the space ${\cal F}$).
\end{proposition}

{\bf Proof:}
\par
{}From the assumed equality: ${\cal D}T^{\tau}_{g}=T^{\tau}_{g}{\cal D}$
it follows easily that ${\cal D}^{-1}T^{\tau}_{g}=T^{\tau}_{g}{\cal D}
^{-1}$. Therefore
$$<\varphi ,T^{\tau}_{g}f> \cong <\eta ,{\cal D}^{-1}T^{\tau}_{g}f> \cong
<\eta ,T^{\tau}_{g}{\cal D}^{-1}f> \cong <\varphi , f>.$$

In the massless case we can consider again equations of type (\ref{spde}),
where now the covariant operator
${\cal D} \in Cov((\R^{N},\tau );(\R^{N},\sigma ))$. The notion
of regularity and the weak solution is defined as in the previous case.

\begin{proposition}
\label{covinter}
Let ${\cal D} \in Cov((\R^{N},\tau );(\R^{N},\sigma ))$ be a regular
operator and let $\eta $ be \par
$\sigma$-covariant random field indexed by
${\cal S}(\R^{D})\otimes \R^{N}$. Then a weak solution of SPDE:
\begin{equation}
\label{spdecovinter}
\tilde{{\cal D}}\varphi = \eta
\end{equation}
given by (\ref{spdesolu}) is $\tau$-covariant random field (providing the
corresponding space ${\cal F}$ is $T^{\tau}_{g}$- invariant).
\end{proposition}

\begin{remark}
Situation such as described in Prop \ref{covinter} occur  for example in
the discussed in \cite{Albcomm} quaternionic representation $\tau = (
\frac{1}{2},\frac{1}{2})$ of $SO(4)$ and the corresponding
Cauchy-Riemann (quaternionic) operator $\partial$.
The corresponding $\sigma = (0,1)$ and the nuclear test function space
${\cal F}$ is defined in Section 3 of \cite{Albcomm}.
In the case of $SO(4)$ it can
be proved that for any irreducible representation $\tau$ the corresponding
set's $Cov(\tau ,\K^{N})$ degenerates to zero-order operators and the only
possibility to produce covariant random fields from SPDE of the type
considered here is to pass to massless case and the choice
${\cal D} \in Cov((\R^{N},\tau );(\R^{N},\sigma ))$. For more details
and new examples in $D=4$ we refer to our forthcoming  paper \cite{Romek1}.
\end{remark}
\begin{remark}
Let $\eta$ be a $\tau$-covariant generalized random field indexed
by ${\cal S}(\R^{D})\otimes \R^{N}$ and let ${\cal D}_{1},...,
{\cal D}_{n},... \in Cov(\tau ,\R^{N})$ be strongly regular
(for simplicity). Let us consider the following cascade of covariant
SPDE's:
\begin{equation}
\tilde{{\cal D}_{n}}\varphi ^{n}=\varphi ^{(n-1)},~~
\tilde{{\cal D}_{1}}\varphi ^{1}=\eta ~~\hbox{for}~n=1,2,3,...
\end{equation}
Then a weak solutions $\varphi ^{n}$ of the cascade $(37)$ (providing
they exist) gives rise to a family $\{ \varphi ^{n}\}$ of $\tau$-covariant
generalized random fields. In particular we have:
\begin{equation}
\Gamma _{\varphi ^{(n)}}(f)=\Gamma _{\eta}({\cal D}_{n}^{-1}...
{\cal D}_{1}^{-1}f)
\end{equation}
\end{remark}
\smallskip
Let ${\cal S}(\R_{+/(-)}\otimes \R^{D-1})=\{ f \in {\cal S}(\R^{D})|
\hbox{supp}f \subset \{x_{0} \geq (\leq) 0,{\bf x} \in \R^{D-1} \}$.\par
Let ${\cal R}:{\cal S}(\R^{D})\otimes \R^{N} \rightarrow
{\cal S}(\R^{D})\otimes \R^{N}$ be a continuous linear mapping such that:
\par
$$(i)~~{\cal R}:
{\cal S}(\R_{+/(-)}^{D})\otimes \R^{N} \rightarrow {\cal S}(\R_{-/(+)}^{D})
\otimes \R^{N}$$
$$(ii)~{\cal R}^{2}=id.$$
A given random field $\eta$ is called ${\cal R}$-reflection positive iff
for all finite sequences $c_{k} \in \C;$ \par
$f^{k} \in {\cal S}(\R_{+}^{D})\otimes \R^{N}$ the following inequality holds:
\begin {equation}
\label{refpos}
\sum_{k,l}c_{k} \overline{c}_{l}\Gamma _{\eta}(f_{k}-{\cal R}f_{l}) \geq 0
\end{equation}

\begin{proposition}
\label{refpositiv}
Let ${\cal D} \in Cov(\tau ,\R^{N})$ be strongly regular and let $\eta$
be ${\cal R}$-reflection positive. Define ${\cal F}_{+/(-)} \equiv
\{ f \in {\cal F}(\R_{+/(-)}^{D})\otimes \R^{N}|f={\cal D}g \hbox{for some}~
g \in {\cal S}(\R_{+/(-)}^{D})\otimes \R^{N}$ and $[{\cal R},{\cal D}]g=0$.
Then the weak solution $\varphi $ of (\ref{spde}) given by (\ref{spdesolu}) is
${\cal R}$-reflection positive in the following sense:
\par
for all finite sequence $c_{k} \in \C$, $f^{k} \in {\cal F}_{+}$:
$$
\sum_{k,l}c_{k} \overline{c}_{l}\Gamma _{\varphi }(f_{k}-{\cal R}f_{l}) \geq 0
$$
\end{proposition}
{\bf Proof}
\par
Let $f^{k}={\cal D}g^{k}$, where $g^{k} \in {\cal S}(\R_{+}\otimes \R^{D-1})
\otimes \R^{N}$. Using ${\cal R}$-reflection positivity (\ref{refpos})
of $\eta$ and
$$\sum_{k,l}c_{k} \overline{c}_{l}\Gamma _{\varphi }(f_{k}-{\cal R}f_{l}) =
\sum_{k,l}c_{k} \overline{c}_{l}\Gamma _{\eta}({\cal D}^{-1}f_{k}-
{\cal D}^{-1}{\cal R}f_{l}) =
\sum_{k,l}c_{k} \overline{c}_{l}\Gamma _{\eta}(g_{k}-{\cal R}g_{l}) \geq 0.
\Box$$

\begin{remark}
Having in mind possible applications of our results to Quantum Field Theory,
examples of ${\cal R}$-reflection positive solutions with suitable reflection
operator ${\cal R}$- has to be produced. We remark that the restricted
reflection positivity demonstrated in Proposition \ref{refpositiv} seems to be
not
sufficiently enough interesting as it leads, in the case where $\eta$ is
taken as physically reflection positive white noise to trivial quantum field
theory models. A detailed discussion of reflection positivity for higher spin
bosonic models of Euclidean Quantum Field Theory together with the proof of
No Go Theorem quoted in the introduction could be find in
\cite{Becker1,Becker2,Romek2}.
\end{remark}

{}From now on we specialize our discussion to the case, when $\eta$ is
$\tau$-covariant white noise with characteristic functional
$\Gamma _{\eta}$ given by $\Gamma _{\eta}=\Gamma _{\eta}^{G} \Gamma _{\eta}
^{P}$ where $\Gamma _{\eta}^{G}$ is given by Gaussian part of
(\ref{Kolmogorovrepresentation}) and
$\Gamma _{\eta}^{P}$ is given by (\ref{ournoise}). We collect some elementary
properties of the weak solution of SPDE (\ref{spde}) with the right hand side
equal
to the white noise as above.
%
%
\par
$<1>$ The weak solution $\varphi$ of regular (\ref{spde}) with ${\cal D} \in
Cov(\tau ,\R^{N})$ has characteristic functional
$\Gamma _{\varphi}$ given by:
\begin{equation}
{\cal F} \ni f \rightarrow \Gamma _{\varphi}^{C}
=\Gamma _{\varphi}^{G}\Gamma _{\varphi}^{P}(f)
\end{equation}
where:
\begin{equation}
\Gamma _{\varphi}^{G}=
e^{-\frac{1}{2}\int <f(x),({\cal D}^{-1})^{T}A{\cal D}^{-1}(x-y)f(x)>}dxdy
\end{equation}
\begin{equation}
\Gamma _{\varphi}^{P}=
e^{\int \int [e^{i<\alpha ,{\cal D}^{-1}\ast f>(x)}-1]}d\nu(\alpha)dx.
\end{equation}
There exists an unique probabilistic Borel cylindric measure
$d\mu _{\cal D}(\varphi)$ on ${\cal F}'(\R^{D}) \equiv$ the weak dual of
${\cal F}$ such that:
\begin{equation}
\Gamma _{\varphi}(f) \equiv \int_{{\cal F}'(\R^{D})} d\mu
_{{\cal D}}(\varphi) e^{i<\varphi ,f>}.
\end{equation}
$<2>$ For any cylinder bounded and of class $C^{1}$ cylindric function
$F \in L^{2}(d\mu _{{\cal D}})$ the following integration by parts
formula holds:
$$
\int_{{\cal F}'(\R^{D})\otimes \R^{N}} <\varphi ,f^{\lambda}>F(\varphi )d\mu
_{{\cal D}}(\varphi) = \int <f^{\lambda}(x), E(({\cal D}^{-1})^{T}A{\cal D}
^{-1}\frac{\delta}{\delta \varphi (x)}F(\varphi))>dx+
$$
\begin{equation}
\int \int f(x) E F
(\varphi +\alpha _{\lambda}\tilde{{\cal D}}^{-1}(\cdot ~-x))\alpha _{\lambda}
d\nu (\alpha)dx.
\end{equation}
where $(f^{\lambda})_{i}=\delta _{i}^{\lambda}f$, $f\in {\cal F}$.
\par
$<3>$ If the Levy measure $d\nu$ has all moments then the field $\varphi$ has
all moments and they are given by the following formula:
\begin{equation}
E\bigl(\,\prod _{i=1}^{n}(\varphi ,f_{i})\bigr)=\sum_{\parbox[t]{50pt}{
$\scriptstyle \Pi_G \cup \Pi_P=J_{n}$ \\
$\scriptstyle \Pi_G \cap \Pi_P = \emptyset $} }
E_P\bigl(\prod _{i \in \Pi_P}(\varphi ,f_{i})\bigr)
E_G\bigl(\prod _{i \in \Pi_G}(\varphi ,f_i)\bigr)
\end{equation}
where:
\begin{equation}
\label{Disturbmom}
E_{P}\bigl(\,\prod _{i=1}^{n}(\varphi ,f_{i})\bigr)
= \sum_{\parbox[t]{60pt}{
$\scriptstyle \Pi _{1}\cup ...
\cup \Pi _{k}=J_{n}$ \\
$\scriptstyle \Pi _{\alpha}\cap \Pi _{\beta}=\emptyset$ \\
$\scriptstyle {\rm for}\; \alpha \neq \beta$} }
\prod _{l=1}^{k} \int ...\int
\prod _{j \in \Pi _{l}}<\alpha_l , {\cal D}^{-1}
f_j (x_l)>\;d\!x_l\;d\!\nu(\alpha_l)
\end{equation}
and
\begin{equation}
E_{G}\bigl(\,\prod_{i=1}^{2n}(\varphi ,f_{i})\bigr)
= \sum_{
  \stackrel{\scriptstyle i_{k} < j_{k}}{
            \scriptstyle k=1,\ldots,n
                                        }
}
\prod _{l=1}^{k} \int \int dxdy
<f_{i_{k}}(x),({\cal D}^{-1})^{T}A{\cal D}^{-1}(x-y)f_{j_{k}}(x)>,
\end{equation}
\begin{equation}
E_{G}\bigl(\,\prod _{i=1}^{2n+1}(\varphi ,f_{i})\bigr)=0.
\end{equation}

In particular the two point moment $S_{\varphi }^{2} \in {\cal F}'^{\otimes 2}$
of $\varphi$ is given by:
\begin{equation}
S_{\varphi }^{2}(f\otimes g)=({\cal D}^{-1})^{T}A{\cal D}^{-1}(f\otimes g)+
\int d\nu (\alpha) \int dx <\alpha ,{\cal D}^{-1}f(x)><\alpha ,{\cal D}^
{-1}g(x)>
\end{equation}
which has the following kernel:
\begin{equation}
S_{\varphi }^{2}(x-y)=({\cal D}^{-1})^{T}A{\cal D}^{-1}(x-y)+
\int d\nu (\alpha) \int dz <\alpha ,{\cal D}^{-1}(z-x)><\alpha ,{\cal D}^
{-1}(z-y)>
\end{equation}
$<4>$ The set $\tilde{{\cal D}}^{-1}\ast C \equiv \{ \sum_{i} \alpha _{i}
\tilde{{\cal D}}^{-1}(\cdot ~-x) |~\hbox{where}~\{ x_{i} \} \in C_{if}(\R^{D})
{}~\hbox{and}~\alpha _{i} \in \hbox{supp}d\nu ~\hbox{for all}~i \}$ is the
carrier set of the Poisson part of the measure $d\mu _{{\cal D}}$ (see
for Remark \ref{carrierProp}).
\par
$<5>$ If the noise is $\tau$-covariant then the random field $\varphi$ is
$\tau$-covariant (providing the test function space ${\cal F}$ is
$T^{\tau}$-invariant).
\par
$<6>$ In the case of strongly regular equation the corresponding solution
is Markovian. The preservation of Markov property under the transformation
$\eta \longrightarrow {\cal D}^{-1}\eta $ with det$\hat{{\cal D}}(ip)\neq 0,
p \in \R^{D}$ follows strightforwardly from paper \cite{Kusuoka}. The case of
nontrivial ker${\cal D}$ is more subtle \cite{Iwata,Piotrek}.

The discussed solutions of SPDE (\ref{spde}) with $\eta$ being Gaussian leads
to
Gaussian (and therefore  not very interesting  from the point of view of
physics) solutions. It is why, we require that the Poisson part of the white
noise  $\eta$ is nonzero in all further applications.
\begin{remark}
Other fundamental properties of the field $\varphi $ like: Markov property,
lattice approximation(s) will be discussed elsewhere (see i.e.
\cite{AlbevI,Iwata,Kusuoka,Nelson1}).
\end{remark}

\section{Laplace-Fourier Transform Properties of the Solutions}
\par
Let us define the following spaces of functions:
$$
{\cal S}_{+}(\R^{Dn}) = \bigl\{ f\in {\cal S}(\R^{Dn}) \ \vert \ f{\rm\ and\
all\ its\ derivatives\ vanish\ unless\ } 0 < x_1^0 < x_2^0 < \ldots <
x_n^0 \bigr\}
$$
$$
{\cal S}_0 (\R^{Dn}) = \bigl\{ f \in {\cal S}(\R^{Dn}) \ \vert \ f{\rm\
and\ all\ its\ derivatives\ vanish\ if\ }x_i= x_j {\rm\ for\ some\ }
1\le i < j\le n \bigr\}
$$
$$
{\cal S}(\R_{+}) = \bigl\{ f\in {\cal S}(\R)\ \vert \ {\rm
supp}\,f\subseteq [0,\infty) \bigr\}, \qquad
{\cal S}(\R_{-}) = \bigl\{ f\in {\cal S}(\R)\ \vert \ {\rm supp}\,
\subseteq (-\infty,0] \bigr\}
.$$
We identify the following spaces:
$${\cal S}(\overline{\R_{+}}) ={\cal S}(\R)/{\cal S}(\R_{-}),~~
{\cal S}(\overline{\R_{+}}^{D}) ={\cal S}(\overline{\R_{+}})\otimes
{\cal S}(\R^{D-1});$$
$${\cal S}(\R^{D};\R^{N}) = \R^{N}\otimes {\cal S}(\R^{D});$$
$${\cal S}(\R^{Dn};(\R^{N})^{\otimes n}) = (\R^{N})^{\otimes n} \otimes
{\cal S}(\R^{Dn});$$
$${\cal S}_{+}(\R^{Dn};(\R^{N})^{\otimes n}) = (\R^{N})^{\otimes n} \otimes
{\cal S}_{+}(\R^{Dn});$$
$${\cal S}_{0}(\R^{Dn};(\R^{N})^{\otimes n}) = (\R^{N})^{\otimes n} \otimes
{\cal S}_{0}(\R^{Dn});$$
$${\cal S}((\overline{\R_{+}}^{D})^{n};(\R^{N})^{\otimes n})
=(\R^{N})^{\otimes n} \otimes
{\cal S}(\overline{\R_{+}}^{Dn}).$$
The following maps will be used:
\begin{equation}
d~:~{\cal S}(\R^{Dn}) \ni f \mapsto f^{d}(x_{1},x_{2}-x_{1},...,x_{n}-x
_{n-1}) \equiv f(x_{1},...,x_{n}).
\end{equation}
The map $d$ is a morphism of ${\cal S}_{+}(\R^{Dn};(\R^{N})^{\otimes n})$
into
${\cal S}((\overline{\R_{+}}^{D})^{n};(\R^{N})^{\otimes n})$. The
Fourier-Laplace transform on ${\cal S}(\overline{\R_{+}}^{Dn})$ is
\begin{equation}
{\cal S}(\overline{\R_{+}}^{Dn}) \ni f_{n} \mapsto f_{n}^{FL}(q_{1},
...,q_{n})=\int e^{-\sum_{k=1}^{n}q_{k}^{0}x_{k}^{0}}e^{i\sum_{k=1}^{n}{\bf q}
_{k}\cdot {\bf x}_{k}}f_{n}(x_{1},...,x_{n}).
\end{equation}
Finally, the map
\begin{equation}
\eta ~:~{\cal S}_{+}(\R^{D(n+1)}) \ni f_{n} \mapsto \eta f_{n} \in
{\cal S}(\overline{\R_{+}}^{Dn})
\end{equation}
is defined as
\begin{equation}
\eta (f_{n})(p_{1},...,p_{n}) \equiv
f_{n}^{d,FL}(p_{1},...,p_{n})\mid _{\{ p_{k}^{0} \geq 0\} }.
\end{equation}
It is well known \cite{oster2} that the map $\eta$ is continuous with
dense range in ${\cal S}(\overline{\R}_{+}^{Dn})$ and trivial kernel.
The notions of $d$, of taking the Fourier-Laplace transform
and of the map $\eta$ naturally extend to the case of
distributions with multiindices.
\begin{definition}
A distribution $F_{n+1} \in {\cal S}_{+}'(\R^{D(n+1)},(\R^{N})^{\otimes (n+1)})$
has the Fourier-Laplace property (the $FL$ property) iff there exists
a distribution ${\cal W}_{n} \in {\cal S}'(\overline{\R}_{+}^{D},(\R^{N})^
{\otimes n})$ such that:
\begin{equation}
F_{n+1}^{d}(x_{0},
...,x_{n}) \equiv \int e^{-\sum_{k=1}^{n}p_{k}^{0}x_{k}^{0}}e^{i\sum_{k=1}^{n}
{\bf p}_{k}\cdot {\bf x}_{k}}{\cal W}_{n}(p_{1},...,p_{n})dp_{1}...dp_{n},
\end{equation}
where the equality $FL$ has to be understand in the sense of distribution
theory sense (see e.g. \cite{Bogoliubov,oster2}).
\end{definition}
There are several necessary and sufficient conditions known for the given
$F_{n} \in {\cal S}_{+}'(\R^{Dn})$ to have $FL$ property
\cite{Bogoliubov,oster2,Simon}. However, all
known for us criterions are hardly to be checked in concrete situations.

Let $\tau$ be a representation of $SO(D)$ in the space $\R^{N}$. We will
say that a tempered distribution $S_{n} \in {\cal S}'(\R^{D};(\R^{N})^
{\otimes n})$ is covariant under the action of $\tau$ ($\tau$-covariant) iff
for each $g \in SO(D)~f_{1},...,f_{n}~\in {\cal S}(\R^{D};\R^{N})$ the
following equality holds:
\begin{equation}
S_{n}(f_{1}\otimes ...\otimes f_{n})=S_{n}(T_{\tau _{g}}f_{1}\otimes...
\otimes T_{\tau _{g}}f_{n})
\end{equation}
where: $(T_{\tau _{g}}f)(x) \equiv \tau _{g} f(g^{-1}x)$.
A distribution $S_{n} \in {\cal S}'(\R^{D};(\R^{N})^{\otimes n})$ is called
symmetric iff
\begin{equation}
S_{n}(f_{1}\otimes ...\otimes f_{n})=S_{n}(f_{\pi (1)}\otimes...
\otimes f_{\pi (n)})
\end{equation}
for any $\pi \in S^{n}(\equiv$ symmetric group of $n$-th element set) and
any $f_{1},...,f_{n} \in {\cal S}(\R^{D};\R^{N})$.
\begin{proposition}
Let $\tau$ be a representation of the group $SO(D)$ in $\R^{N}$. If
$\sigma _{n} \! \in \! {\cal S}'(\R^{D};(\R^{N})^{\otimes n})$ is symmetric
covariant
under the action of $\tau$ and $\sigma _{n}\mid _{\{ y^{0}_{k} \geq 0\} }$
has $FL$
property then there exists a unique tempered distribution ${\cal W}_{n}
\in {\cal S}'(\R^{D};(\R^{N})^{\otimes n})$ such that:
\par
(1) ${\cal W}_{n}^{F}$ is supported in the product of forward light
cones $V^{+}\equiv \{p \in M^{D}|p\cdot p \geq 0;p^{0} \geq 0 \}$, i.e.:
$$supp{\cal W}_{n} \subseteq (V^{+})^{\times n}.$$
(2) ${\cal W}_{n}$ is covariant under the representation $\tau ^{M}$ of
$SO(D-1,1)$, i.e.:
\begin{equation}
{\cal W}_{n}(f_{1}\otimes ...\otimes f_{n})={\cal W}_{n}
(T_{\tau _{g}^{M}}f_{1}\otimes...\otimes T_{\tau _{g}^{M}}f_{n})
\end{equation}
for any $g \in L_{+}^{\uparrow}(D);f_{1},...,f_{n}\in {\cal S}(\R^{D};R^{n})$
and where $\tau ^{M}$ is the analytic continuation of $\tau$ into the
representation of $SO(D-1,1)$ via the "Weyl unitary trick".
\par
(3) ${\cal W}_{n}$ is local, which means, that the inverse Fourier transform of
${\cal W}_{n}(x_{1},...,x_{n})$ has the property that
if some $x_{i},x_{i+1}$ are  such that $(x_{i}-x_{i+1})^{2}<0$ then
\begin{equation}
{\cal W}_{n}(x_{1},...,x_{i},x_{i+1},...,x_{n})=
{\cal W}_{n}(x_{1},...,x_{i+1},x_{i},...,x_{n})
\end{equation}
(4)\begin{equation}
S_{n+1}^{d}(x_{0},...,x_{n})
 \equiv \int e^{-\sum_{k=1}^{n}p_{k}^{0}x_{k}^{0}}e^{i\sum_{k=1}^{n}
{\bf p}_{k}\cdot {\bf x}_{k}}{\cal W}_{n}^{F}(p_{1},...,p_{n})
\prod _{i=1}^{n}dp_{i}
\end{equation}
for $x_{1}^{0}\leq ... \leq x_{n}^{0}$.
\end{proposition}
{\bf Proof:}
\par
{}From the $FL$ property of $S_{n}$ it follows that there exists ${\cal W}_{n}
\in {\cal S}'(\R^{D};(\R^{N})^{\otimes n})$ such that ${\cal W}_{n}$ is
supported on positive  energies, i.e. on the set $\{ (p_{1},...,p_{n})|~
p_{i}^{0} \geq 0 ~$ for all $i=1,...,n \}$ and such that $(4)$ holds. But
the covariant under the action of the Lorentz group distribution must be
supported in the orbit of Lorentz group
\cite{Bogoliubov,Jost,Streater} and thus we conclude that
${\cal W}_{n}$ must be supported on $(V^{+})^{\times n}$. The locality of
${\cal W}_{n}$ follows from the symmetry property of $S_{n}$ (see e.g.
\cite{Jost}).
The uniqueness of ${\cal W}_{n}$ follows from the fact that the kernel of the
Laplace-Fourier transform consists only from vector $0$. $\Box$

The difference variables moments $\sigma _{n}$ of the random fields $A$
constructed in section 2 are defined as:
\begin{equation}
\sigma _{n}^{A}(\xi _{1},...,\xi _{n}) \equiv S_{n+1}^{A}(x_{1},...,
x_{n+1})
\end{equation}
where $\xi _{i} \equiv x_{i+1}-x_{i}$ for $i=1,...,n$ 
Now we are ready to formulate the main result of this paper.

\begin{theorem}
Let $\tau$ be a real representation of $SO(D)$ in $\R^{N}$, ${\cal D} \in
Cov(\tau,\R^{N})$ with admissible spectrum and let $A$ be a solution of
$$\tilde{{\cal D}}A~=~\eta $$
where $\eta$ is the $T_{\tau }$-invariant Poisson noise. Then the
difference variables moments $\sigma _{n}^{A}(x_{1},...,x_{n})$ have
Fourier-Laplace property.
\end{theorem}
The proof of this theorem will be divided into three main steps.
\begin{proposition}
Let ${\cal D} \in Cov(\tau ,\R^{N})$ has an admissible mass spectrum with
strictly positive masses.
Then the Green function $G_{A}=({\cal D}_{A})^{-1\ast}$ of ${\cal D}$ has
Fourier-Laplace property.
\end{proposition}
\begin{lemma}
Let $A$, $G_{A}$ be as in Proposition 3-3.
\begin{equation}
S'(|x_{1}-x_{2}|)=\int dx{\cal D}_{A}^
{-1}(x-x_{1}){\cal D}_{A}^{-1}(x-x_{2})
\end{equation}
has the Fourier-Laplace property.
\end{lemma}
\begin{lemma}
Let $A$, $G_{A}$ be as in Proposition 3-4. Then for any $k=1,2$ the
distribution
$S'^{d}_{k}$, where
\begin{equation}
S_{k}^{,d}(x_{1},...,x_{k}) \equiv \int dx{\cal D}_{A}^{-1}(x-x_{1})...
{\cal D}_{A}^{-1}(x-x_{k})
\end{equation}
has the Fourier-Laplace property.
\end{lemma}
The separation of the proof into Lemma 3-5 and Lemma 3-6 is made for reader
convenience only.
Having proven Proposition 3-4, Lemma 3-6 the proof of theorem 3-3 follows by
noting the fact that the Fourier-Laplace property is stable under taking tensor
product and use of formula (\ref{Disturbmom}).
The case in which some of the masses are
equal to zero is easily covered by using the continuity of Laplace-Fourier
transform in the space of distributions and a easily controlled limiting
procedure based on an introducing virtual nonzero masses in the corresponding
formulae and then putting them to zeros. Although the covariance might be
broken by introducing virtual masses it can be restored in the limit.
\par
{\bf Proof} of Proposition 3-4.
\par

The typical matrix element $G_{A}^{\alpha \beta}$ of $G_{A}$ has the form
\begin{equation}
G_{A}^{\alpha \beta}(p)=\frac{Q^{\alpha \beta }(p)}{\prod _{i=1}^{n}
(p_{0}^{2}+{\bf p}^{2}+m_{i}^{2})}
\end{equation}
where $Q^{\alpha \beta}$ are polynomials in variables $p$ of degree lower or
equal
$N-2$ and all $m_{i}>0$ due to assumption made on ${\cal D}$;
\begin{equation}
\label{fractions}
\frac{Q^{\alpha \beta }(p)}{\prod _{i=1}^{n}
(p_{0}^{2}+{\bf p}^{2}+m_{i}^{2})}=\sum_{i=1}^{n}\frac{Q^{\alpha \beta }_{i}
(p_{0},{\bf p})}{(p_{0}^{2}+{\bf p}^{2}+m_{i}^{2})}
\end{equation}
where $Q^{\alpha \beta}_{i}(p_{0},{\bf p}) \equiv A^{\alpha \beta ,i}({\bf p})
p_{0}+B^{\alpha \beta ,i}({\bf p})$ where $A^{\alpha \beta ,i}({\bf p}),
B^{\alpha \beta ,i}({\bf p})$ are bounded (on $\R$) rational functions in
variable ${\bf p}$.

As it is well know the distribution
\begin{equation}
S_{\Box ,i}^{2}(x)=\int \frac{1}{p_{0}^{2}+{\bf p}^{2}+m_{i}^{2}}
e^{-ipx}dp
\end{equation}
has the Fourier-Laplace property with the underlying distribution ${\cal W}^{0}
_{i}$ given by ${\cal W}^{0}_{i}(p_{0},{\bf p})=\epsilon (p_{0})\delta (p_{0}
^{2}-{\bf p}^{2}-m^{2})$, where $\epsilon (p_{0})=1$ if $p_{0} \geq 0$ and $0$
otherwise. The inverse Fourier transform of a typical term appearing in
(\ref{fractions}) is
given (for $x^{0} \geq 0$):
$$(A^{\alpha \beta}_{i}(i{\bf \nabla})i\frac{\partial}{\partial x^{0}}+
B^{\alpha \beta}_{i})S^{2}_{\Box ,i}(x^{0},{\bf x})=
(A^{\alpha \beta}_{i}(i{\bf \nabla})i\frac{\partial}{\partial x^{0}}+
B^{\alpha \beta}_{i})(\! \int_{0}^{\infty} \!
\int \! e^{-p_{0}x_{0}}e^{i{\bf p}\cdot {\bf x}}
\delta (p_{0}^{2}-{\bf p}^{2}-m^{2})dp_{0}d{\bf p})$$
$$
=B^{\alpha \beta}_{i}
\! \int_{0}^{\infty}\! \int \! e^{-p_{0}x_{0}}e^{i{\bf p}\cdot {\bf x}}
\delta (p_{0}^{2}-{\bf p}^{2}-m^{2})dp_{0}d{\bf p}+$$
\begin{equation}
+\! \int_{0}^{\infty} \! \int \! e^{-p_{0}x_{0}}e^{i{\bf p}\cdot {\bf x}}
\{ -ip_{0}A^{\alpha \beta}_{i}(i{\bf p})\delta (p_{0}^{2}-{\bf p}^{2}-m^{2})\}
dp_{0}d{\bf p}
\end{equation}
and this shows that the inverse Fourier transform of each term in
(\ref{fractions}) is the
Fourier-Laplace transform with underlying distribution for $i$-th term:
\begin{equation}
{\cal W}^{\alpha \beta}_{i}(p_{0},{\bf p})=\{ B^{\alpha \beta}_{i}
-ip_{0}A^{\alpha \beta}_{i}(i{\bf p})\} \delta (p_{0}^{2}-{\bf p}^{2}-m^{2})
\epsilon (p_{0})
\end{equation}
and ${\cal W}^{\alpha \beta}(p_{0},{\bf p}) \equiv
\sum_{i} {\cal W}_{i}^{\alpha \beta}(p_{0},{\bf p})$. $\Box$
\par
{\bf Proof} of Lemma 3-5
\par

We shall proceed very close to the proof of Thm 4.21 in \cite{Albcomm}.
Firstly, we use the following identity:
\begin{equation}
\int_{-\infty}^{+\infty}e^{-\zeta _{1}|t-t_{1}|}e^{-\zeta _{2}|t-t_{2}|}dt=
\frac{1}{\zeta _{1}+\zeta _{2}}e^{-\zeta _{2}(t_{2}-t_{1})}+
\frac{1}{\zeta _{1}+\zeta _{2}}e^{-\zeta _{1}(t_{2}-t_{1})}+
(t_{2}-t_{1})\cdot \int_{0}^{1}e^{-(\zeta _{1}s+(1-s)\zeta _{2})(t_{2}-t_{1})}
\end{equation}
which is valid for any $t_{1},t_{2} \in \R;~\zeta _{1},\zeta _{2} \in \C$
such that: $t_{2}-t_{1}>0$; $\Re \zeta _{1}>0$; $\Re \zeta _{2}>0$.
Secondly we note that:
$$\int_{0}^{\infty}e^{-p_{0}|x_{0}|}e^{-i{\bf p}\cdot {\bf x}}{\cal W}^{
\alpha \beta }_{i}(p_{0},{\bf p})dp_{0}d{\bf p}=-\frac{i}{(2\pi )^{D-1}}
\int_{\R^{D-1}}e^{-\sqrt{{\bf p}^{2}+m_{i}^{2}}|x_{0}|}e^{-i{\bf p}\cdot
{\bf x}}A^{\alpha \beta}_{i}({\bf p})d^{D-1}{\bf p}$$
\begin{equation}
+B^{\alpha \beta}_{i}\int_{\R^{D-1}}
\frac{e^{-\sqrt{{\bf p}^{2}+m_{i}^{2}}|x_{0}|}}{\sqrt{{\bf p}^{2}+m_{i}^{2}}}
e^{-i{\bf p}\cdot {\bf x}}A^{\alpha \beta}_{i}({\bf p})d^{D-1}{\bf p}.
\end{equation}
Now, we can write down:
\begin{equation}
\Gamma (y_{1}-y_{2}) \equiv \int dx^{0}d{\bf x} {\cal D}^{-1}(x-y_{1})
{\cal D}^{-1}(x-y_{2}) =
\end{equation}
$$
\int dx^{0}d{\bf x} {\cal D}^{-1}(|x^{0}-y^{0}_{1}|,|{\bf x}-{\bf y}_{1}|)
{\cal D}^{-1}(|x^{0}-y^{0}_{2}|,|{\bf x}-{\bf y}_{2}|)=
$$
$$
\int_{-\infty}^{+\infty}\! dx^{0} \!
\int \! d{\bf x} \! \int \! dp_{1} \! \int \! dp_{2} \! e^{-p^{0}_{1}
|x^{0}-y^{0}_{1}|}e^{-p_{2}^{0}|x^{0}-y^{0}_{2}|}
e^{-i{\bf p}_{1}\cdot ({\bf x}-{\bf y}_{1})}e^{-i{\bf p}_{2}\cdot ({\bf x}-
{\bf y}_{2})}{\cal W}_{G}(p_{1}^{0},{\bf p}_{1}){\cal W}_{G}(p_{2}^{0},
{\bf p}_{2})
$$
$$
=\int \! d{\bf x} \! \int \! dp_{1} \! \int \! dp_{2} \!
\{ \frac{1}{p_{1}^{0}+
p_{2}^{0}}(e^{-p^{0}_{1}(y_{2}^{0}-y_{1}^{0})}+
e^{-p^{0}_{2}(y_{2}^{0}-y_{1}^{0})})+
$$
$$
+(y_{2}^{0}-y_{1}^{0})\int_{0}^{1} \! ds \!
e^{-(p_{1}^{0}s+(1-s)p_{2}^{0})(y_{2}^{0}-y_{1}^{0})}
e^{-i{\bf p}_{1}\cdot ({\bf x}-{\bf y}_{1})}e^{-i{\bf p}_{2}\cdot ({\bf x}-
{\bf y}_{2})}{\cal W}_{G}(p_{1}^{0},{\bf p}_{1}){\cal W}_{G}(p_{2}^{0},
{\bf p}_{2})\}$$
$$ \equiv 2\Gamma _{1}(y_{1}-y_{2})+\Gamma _{2}(y_{1}-y_{2}).$$
Defining the following functions
$$
\Pi _{(1)}^{\alpha \beta i,\alpha '\beta 'i^{\ast}}=
\int d{\bf p}e^{-i{\bf p}\cdot ({\bf y}_{2}-{\bf y}_{1})}
\frac{\sqrt{{\bf p}^{2}+m_{i}^{2}}\sqrt{{\bf p}^{2}+m_{i^{\ast}}^{2}}}
{\sqrt{{\bf p}^{2}+m_{i}^{2}}+\sqrt{{\bf p}^{2}+m_{i^{\ast}}^{2}}}
\{ e^{-\sqrt{{\bf p}^{2}+m_{i}^{2}}(y_{2}^{0}-y_{1}^{0})}+ $$
\begin{equation}
\label{ser1}
+ e^{-\sqrt{{\bf p}^{2}+m_{i^{\ast}}^{2}}(y_{2}^{0}-y_{1}^{0})} \}
A^{\alpha \beta}_{i}(-i{\bf p})A^{\alpha '\beta '}_{i^{\ast}}(i{\bf p})
\end{equation}
\begin{equation}
\Pi _{(2)}^{\alpha \beta i,\alpha '\beta 'i^{\ast}}=
\! \int \! d{\bf p} e^{-i{\bf p}\cdot ({\bf y}_{2}-{\bf y}_{1})}
\frac{B^{\alpha \beta}_{i}B^{\alpha '\beta '}_{i^{\ast}}}
{\sqrt{{\bf p}^{2}+m_{i}^{2}}+\sqrt{{\bf p}^{2}+m_{i^{\ast}}^{2}}}
\{ e^{-\sqrt{{\bf p}^{2}+m_{i}^{2}}(y_{2}^{0}-y_{1}^{0})}+
e^{-\sqrt{{\bf p}^{2}+m_{i^{\ast}}^{2}}(y_{2}^{0}-y_{1}^{0})}\}
\end{equation}
\begin{equation}
\Pi _{(3)}^{\alpha \beta i,\alpha '\beta 'i^{\ast}}=
\int d{\bf p}e^{-i{\bf p}\cdot ({\bf y}_{2}-{\bf y}_{1})}
\frac{B^{\alpha \beta }_{i}(-i)\sqrt{{\bf p}^{2}+m_{i^{\ast}}^{2}}}
{\sqrt{{\bf p}^{2}+m_{i}^{2}}+\sqrt{{\bf p}^{2}+m_{i^{\ast}}^{2}}}
A^{\alpha '\beta '}_{i^{\ast}}(i{\bf p})
\end{equation}
$$
\Pi _{(4)}^{\alpha \beta i,\alpha '\beta 'i^{\ast}}=
\int d{\bf p}e^{-i{\bf p}\cdot ({\bf y}_{2}-{\bf y}_{1})}
B^{\alpha '\beta '}_{i^{\ast}}(-i)\sqrt{{\bf p}^{2}+m_{i}^{2}}
\{ e^{-\sqrt{{\bf p}^{2}+m_{i}^{2}}(y_{2}^{0}-y_{1}^{0})}+$$
\begin{equation}
+e^{-\sqrt{{\bf p}^{2}+m_{i^{\ast}}^{2}}(y_{2}^{0}-y_{1}^{0})} \}
A^{\alpha \beta}_{i}(i{\bf p})
\end{equation}
$$
\Gamma _{(1)}^{\alpha \beta i,\alpha '\beta 'i^{\ast}}=
\int d{\bf p} e^{-i{\bf p}\cdot ({\bf x}-{\bf y})} \int_{0}^{1}ds
(y_{2}^{0}-y_{1}^{0})\sqrt{{\bf p}^{2}+m_{i}^{2}}
\sqrt{{\bf p}^{2}+m_{i^{\ast}}^{2}} \times $$
\begin{equation}
A^{\alpha \beta}_{i}(-i{\bf p}_{1})
A^{\alpha \beta}_{i^{\ast}}(i{\bf p}_{2})e^{-\sqrt{{\bf p}^{2}+m_{i}^{2}}
(y_{2}^{0}-y_{1}^{0})s}e^{-\sqrt{{\bf p}^{2}+m_{i}^{2}}(y_{2}^{0}-y_{1}^{0})
(1-s)}
\end{equation}
\begin{equation}
\Gamma _{(2)}^{\alpha \beta i,\alpha '\beta 'i^{\ast}}=
B^{\alpha \beta}_{i}B^{\alpha \beta}_{i^{\ast}}(y_{2}^{0}-y_{1}^{0})
\int d{\bf p} e^{-i{\bf p}\cdot ({\bf x}-{\bf y})} \int_{0}^{1}ds
e^{-\sqrt{{\bf p}^{2}+m_{i}^{2}}
(y_{2}^{0}-y_{1}^{0})s}e^{-\sqrt{{\bf p}^{2}+m_{i}^{2}}(y_{2}^{0}-y_{1}^{0})
(1-s)}
\end{equation}
$$
\Gamma _{(3)}^{\alpha \beta i,\alpha '\beta 'i^{\ast}}=
(y_{2}^{0}-y_{1}^{0})
\int d{\bf p} e^{-i{\bf p}\cdot ({\bf y}_{2}-{\bf y}_{1})}(-i)
\sqrt{{\bf p}^{2}+m_{i}^{2}}A^{\alpha '\beta '}_{i^{\ast}}(i{\bf p})
\int_{0}^{1}ds \times $$
\begin{equation}
e^{-\sqrt{{\bf p}^{2}+m_{i}^{2}}(y_{2}^{0}-y_{1}^{0})s}
e^{-\sqrt{{\bf p}^{2}+m_{i}^{2}}(y_{2}^{0}-y_{1}^{0})(1-s)}
\end{equation}
$$
\Gamma _{(4)}^{\alpha \beta i,\alpha '\beta 'i^{\ast}}=
(y_{2}^{0}-y_{1}^{0})
\int d{\bf p} e^{-i{\bf p}\cdot ({\bf y}_{2}-{\bf y}_{1})}(-i)
\sqrt{{\bf p}^{2}+m_{i^{\ast}}^{2}}A^{\alpha '\beta '}_{i}(i{\bf p})
\int_{0}^{1}ds\times
$$
\begin{equation}
\label{ser2}
e^{-\sqrt{{\bf p}^{2}+m_{i}^{2}}(y_{2}^{0}-y_{1}^{0})s}
e^{-\sqrt{{\bf p}^{2}+m_{i}^{2}}(y_{2}^{0}-y_{1}^{0})(1-s)}.
\end{equation}
We obtain after a bit of calculations that:
\begin{equation}
S_{,2}^{\alpha \beta ,\alpha '\beta '}(y_{2}-y_{1})=\int dx {\cal D}^{-1}_
{\alpha \beta }(x-y_{1}){\cal D}^{-1}_{\alpha '\beta '}(x-y_{2}) \equiv
\sum_{\delta =1}^{4}(\sum_{i,i^{\ast}}\Pi _{\delta}
^{\alpha \beta i,\alpha '\beta 'i^{\ast}} + \sum_{i,i^{\ast}} \Gamma _{\delta}
^{\alpha \beta i,\alpha '\beta 'i^{\ast}}).
\end{equation}
{}From the explicite formulae (\ref{ser1})-(\ref{ser2})
it follows that all the functions possess an analytic continuation.
\par
{\bf Proof} of lemma 3.6:
\par
The following (see eqs. 4.16 in \cite{Albcomm})
$$\int \prod_{i=1}^{n}e^{-\zeta _{i}|t-t_{i}|}dt=\frac{1}{\zeta _{1}+...+
\zeta _{n}}e^{-\zeta _{2}(t_{2}-t_{1})}...e^{-\zeta _{n}(t_{n}-t_{1})}+$$
$$
+\sum_{j=1}^{n-1}\prod_{i=1}^{j-1}e^{-\zeta _{i}(t_{j}-t_{i})}(t_{j+1}-t_{i})
\prod_{i=j+2}^{n}e^{-\zeta _{i}(t_{i}-t_{j+1})}\int
e^{-[(\zeta _{1}+...+\zeta _{j})s+(\zeta _{j+1}+...+\zeta _{n})(1-s)(t_{j+1}-
t_{j})]}ds +$$
\begin{equation}
+ \frac{1}{\zeta _{1}+...+\zeta _{n}}e^{-\zeta _{1}(t_{n}-t_{1})}
...e^{-\zeta _{n-1}(t_{n}-t_{n-1})}
\end{equation}
is valid for any $t_{1}<t_{2}<...<t_{n}$ and complex numbers $\zeta _{i}$ such
that $\Re \zeta _{i} >0$ for all $i$ and the decomposition
(\ref{fractions}) is used to
derive the following representation of $S_{k}^{,d}$:
$$
S_{k}^{,d~\alpha _{1}\beta _{1}...\alpha _{k}\beta _{k}}(x_{1},...,x_{k})
\equiv \int dx {\cal D}_{A}^{-1~~\alpha _{1}\beta _{1}}(x-x_{1})...
{\cal D}_{A}^{-1~~\alpha _{k}\beta _{k}}(x-x_{k})=
$$
$$
=\sum_{\delta _{1}=1}^{n}...\sum_{\delta _{n}=1}^{n}\int dp_{1}^{0}d{\bf p}
_{1}...dp_{k}^{0}d{\bf p}_{k} \prod_{j=1}^{k}{\cal W}_{\delta _{k}}^{\alpha
_{k}\beta _{k}}(p_{j}^{0},{\bf p}_{j}) \int d{\bf x} \prod_{j=1}^{k}
e^{-i({\bf x}-{\bf x}_{j})\cdot {\bf p}_{j}} \int dx^{0}
\prod_{j=1}^{k}e^{-p_{j}^{0}|x^{0}-x_{j}^{0}|}=
$$
$$
\sum_{\delta _{1}=1}^{n} \ldots \sum_{\delta _{n}=1}^{n} \! \int \! dp_{1}^{0}
d{\bf p}_{1} \ldots
dp_{k}^{0}d{\bf p}_{k} \! \int \! d{\bf x} \! \prod_{j=1}^{k} \!
e^{-i({\bf x}-{\bf x}_{j})\cdot {\bf p}_{j}}
\! \prod_{j=1}^{k} \! {\cal W}_{\delta _{k}}^{\alpha _{k}
\beta _{k}}(p_{j}^{0},{\bf p}_{j})
\{ \frac{1}{p_{1}^{0}+...+p_{k}^{0}}
\prod_{j=1}^{k}e^{-p_{j}^{0}(x_{j}^{0}-x_{1}^{0})}$$
$$+\sum_{j=1}^{n-1}\prod_{i=1}^{j-1}e^{-p_{i}^{0}(x_{j}^{0}-x_{i}^{0})}
(x_{j+1}^{0}-x_{i}^{0}) \! \prod_{i=j+2}^{n} \! e^{-p_{i}^{0}
(x_{i}^{0}-x_{j+1}^{0})} \! \int \! ds \! e^{-[(p_{1}^{0}+...+p_{j}^{0})s+
(p_{j+1}^{0}+...+p_{k}^{0})(1-s)(x_{j+1}^{0}-x_{j}^{0})]} +$$
\begin{equation}
+\frac{1}{p_{1}^{0}+...+p_{k}^{0}}e^{-p_{1}^{0}(x_{k}^{0}-x_{1}^{0})}
...e^{-p_{k-1}^{0}(x_{k}^{0}-x_{k-1}^{0})} \}.
\end{equation}
Similarly, as in the proof of lemma 3.5, when using the explicite expressions
for $\{ {\cal W}_{\delta}^{\alpha \beta}(p_{0},{\bf p})\}$ given in the proof
of Proposition 3.4 one can see that the functions $S_{k}'^{d~\alpha _{1}
\beta _{1}...\alpha _{k}\beta _{k}}$ are given by sums, each term of whose
is manifestly given by the Fourier-Laplace transforms of some tempered
distribution supported on positive energies.
\begin{remark}
Let $\{ {\cal W}_{n} \}$ be the obtained set of $\tau_{M}$-covariant, local and
fulfilling the weak form of the spectral axiom Wightman distributions. Then
using a version of GNS construction one could construct: an inner product space
${\cal H}^{ph}$ with the inner product $<\cdot,\cdot>_{{\cal H}^{ph}}$, a
linear
weakly continuous map:
$$A_{q} : {\cal S}(\R^{D})\otimes \R^{N} \longrightarrow \ell ({\cal H}^{ph})$$
where $\ell ({\cal H}^{ph})\equiv$ the set of linear (not necessarily bounded)
operators acting on ${\cal H}^{ph}$, nonunitary and unbounded representation
$U^{M}_{\tau _{M}}$ of ${\cal P}_{+}^{\uparrow}(D)$ in ${\cal H}^{ph}$
under which the quantum field operator $A_{q}$ transforms covariantly, and a
cyclic with respect to the action of $A_{q}(f)$ and invariant with respect to
$U^{M}_{\tau _{M}}$ vector $\Omega$ playing role of physical vaccum.
\end{remark}

\section{Examples in $D=3$}

The complete description of the set of all covariant operators
${\cal D} \in Cov(\tau)$, where $\tau$ is any finite dimensional representation
of the group $SO(3)$ or $SO(1,3)$ is given in the monographs
\cite{GelfandShapiro,Naimark} (see also \cite{Ljubarski}).
To illustrate our general theory developed in the previous
paragraphs we focus attention on the lowest-dimensional real representations
$D_{0}\oplus D_{1}$, $D_{1}\oplus D_{1}$ and $D_{\frac{1}{2}}\oplus D_{
\frac{1}{2}}$ of the group $SO(3)$. The much more interesting case of $D=4$
shall be analysed in a greater details in our forthcoming paper \cite{Romek1}.
The presented below examples do not exhaust all the possibilities. The
point is, that we have used rather special realification procedure in order
to brought the complex description of the sets $Cov(\tau)$ given in
\cite{GelfandShapiro,Naimark,Ljubarski}
into the manifestly real form. Our realification is achieved
by certain similarity transformation, fixed by the choice of a
realification matrix
$E_{\tau}$. [Different choices of the realification procedure may lead to
a different (i.e. not connected by the similarity transformation) families of
covariant operators].
\subsection{$D_0\oplus D_1$: Higgs-like Models}

This class of models described a doublet of fields $\varphi =(\varphi _{0},
{\bf A})$ where $\varphi _{0}$ is the scalar field and ${\bf A}$ is the vector
field, coupled by noise throughout the corresponding covariant SPDE of the form
(\ref{spde}).

The realification matrix $E_{(0,1)}$ is chosen to be:
\begin{equation}
              E_{(0,1)}= \frac{1}{\sqrt{2}}
              \left( \begin{array}{cccc}
                      \sqrt{2} & 0 & 0 & 0  \\
                             0 & i & 0 & -i \\
                             0 & 1 & 0 & 1  \\
                             0 & 0 & i\sqrt{2} & 0
                                \end{array} \right)
\end{equation}
The real form of the corresponding covariant operators ${\cal D}_{(0,1)}
\in Cov((D_{0}\oplus D_{1})^{R})$ with respect to $SO(3)$ with the mass
term $M=m_{0}{\bf 1}_{0} \oplus m_{1}{\bf 1}_{3}$ is given by:
\begin{equation}
\label{allcov}
             \hat{\cal D} _{(0,1)}= \left( \begin{array}{cccc}
               m_{0} & aip_{0} & aip_{1}  &  aip_{2} \\
             bip_{0} &   m_{1} & -cip_{2} &  cip_{1} \\
             bip_{1} & cip_{2} &    m_{1} & -cip_{0} \\
             bip_{2} & -cip_{1} & cip_{0} &    m_{1}
                                       \end{array} \right)
\end{equation}
with $a,b,c \in \R$ with $\hbox{det}\hat{\cal D}
_{(0,1)}(ip)=(-c^{2}p^{2}+m_{1}
^{2})(abp^{2}+m_{0}m_{1})$

To obtain admissible mass spectrum we need to put either $c=0$ or $m_{1}=0$,
therefore
resignating from the ellipticity of $\hat{{\cal D}}$. An admissible covariant
operators are obtained iff $c=0$ if $m_{1} \neq 0$ or $m_{1}=0$.
We add that the operator is covariant with respect to $O(3)$ iff $c=0$.
The Green function is given by:
\begin{equation}
\hat{{\cal D}}_{(0,1)}^{-1}(ip)= \frac{1}{abp^{2}+m_{0}m_{1}}\left(
       \begin{array}{c|c}
        m_{1} & \begin{array}{ccc} -aip_{0} & -aip_{1} & -aip_{2}
                  \end{array}
                                    \\ \hline
     \begin{array}{c} -bip_{0} \\  -bip_{1} \\ -bip_{2}
                  \end{array}  & G_{\mu \nu}(p)
       \end{array} \right)
\end{equation}
where:
$$G_{\mu \nu}(p)=\frac{1}{-c^{2}p^{2}+m_{1}^{2}}\{ (abp^{2}+m_{0}m_{1})
(m_{1}\delta _{\mu \nu} + ci\varepsilon _{\mu \nu \lambda} p_{\lambda})
-p_{\mu}p_{\nu}(abm_{1}+c^{2}m_{0}) \} $$
for $\mu ,\nu \in \{ 0,1,2 \}$. The corresponding two-point function
(more precisely the contribution coming from the Poisson piece of noise
and not integrated with L\'evy measure $\nu$, see eq. $2.50$):
\begin{equation}
\hat{S}_{(0,1)}^{(2)}(p,\alpha)=
(\hat{S}^{(2)}_{kl}(p,\alpha))= \left(
                  \begin{array}{c|c}
   \hat{S}^{(2)}_{33}(p,\alpha) &
                    \begin{array}{ccc}
\hat{S}^{(2)}_{30}(p,\alpha) & \hat{S}^{(2)}_{31}(p,\alpha) &
              \hat{S}^{(2)}_{32}(p,\alpha)
                    \end{array}
                                     \\ \hline
            \begin{array}{c}
 \hat{S}^{(2)}_{03}(p,\alpha) \\ \hat{S}^{(2)}_{13}(p,\alpha) \\
\hat{S}^{(2)}_{23}(p,\alpha)
            \end{array}
            & \hat{S}^{(2)}_{\mu \nu}(p,\alpha)
                  \end{array} \right)
\end{equation}
where
$$\hat{S}^{(2)}_{33}(p,\alpha)=
|m_{1}\alpha _{3}-ib\alpha _{\mu}p_{\mu}|^{2}/(abp^{2}+m_{0}m_{1})^{2},
$$
$$
\hat{S}^{(2)}_{3 \mu}(p,\alpha)=\hat{S}^{(2)}_{\mu 3}(p,\alpha)=
(m_{1}\alpha_{3}-ib\alpha_{\lambda}p_{\lambda})[m_{1}(abp^{2}+m_{0}m_{1})
\alpha_{\mu}+ia(-c^{2}p^{2}+m_{1}^{2})\alpha_{3}p_{\mu}-
$$
$$-(abm_{1}+c^{2}m_{0})
\alpha _{\lambda}p_{\lambda}p_{\mu}]/(-c^{2}p^{2}+m_{1}^{2})
(abp^{2}+m_{0}m_{1})^{2},
$$
$$
\hat{S}^{(2)}_{\mu \nu}(p,\alpha)=
[a^{2}\alpha _{3}^{2}+(\alpha _{\lambda}p_{\lambda})^{2}(abm_{1}+
c^{2}m_{0})^{2}
/(-c^{2}p^{2}+m_{1}^{2})^{2}]p_{\mu}p_{\nu}/(abp^{2}+m_{0}m_{1})^{2}+
$$
$$
m_{1}^{2}\alpha_{\mu}\alpha_{\nu}/(-c^{2}p^{2}+m_{1}^{2})^{2}-m_{1}
(abm_{1}+c^{2}m_{0})\alpha _{\lambda}p_{\lambda}
(p_{\mu}\alpha _{\nu}+p_{\nu}\alpha_{\mu})/(-c^{2}p^{2}+m_{1}^{2})^{2}
(abp^{2}+m_{0}m_{1})^{2}-
$$
$$-iam_{1}\alpha_{3}
(p_{\mu}\alpha _{\nu}-p_{\nu}\alpha_{\mu})/(-c^{2}p^{2}+m_{1}^{2})
(abp^{2}+m_{0}m_{1})
$$
for $\mu ,\nu \in \{ 0,1,2 \}$. \par
We use above the notation for the variable
of L\'evy measure: $\alpha \equiv (\alpha _{3}, \alpha _{0}, \alpha _{1},
\alpha _{2})$.\par
{\bf Remarks} \par
The representation $D_{0}\oplus D_{1}$ is also of quaternionic type.
Choosing $m_{0}^{2}+m_{1}^{2}=0$ and $a=-1,b=1,c=1$ (respectively
$a=-1,b=1,c=-1$) in $(2)$ we obtain the purely quaternionic description of the
corresponding Clifford algebra of $\R^{3}$ Dirac operators. More explicitly
let denote $C(\R^{3})$ the corresponding to $\R^{3}$ Clifford algebra and by
$\Lambda (\R^{3})$ the external algebra of $\R^{3}$. Let us denote by
$C(\R^{3})=C_{+}(\R^{3})\oplus C_{-}(\R^{3})$ (respectively
$\Lambda (\R^{3})=\Lambda _{+}(\R^{3})\oplus \Lambda _{-}(\R^{3})$)
canonical decompositions of $C(\R^{3})$ (resp. of $\Lambda (\R^{3})$)
gradation. Let $\HH$ stands for the noncommutative field of quaternions
with the base $\{ 1,{\bf i},{\bf j},{\bf k} \}$. Noting that
$C(\R^{3}) \cong \HH \oplus \HH $ canonical and $\Lambda (\R^{3}) \cong
\HH \oplus \HH$ and using  two (non-equivalent) representations of $\HH$
on $\HH$ given by left (resp. right) multiplication we obtain the following
explicite expressions for the corresponding left (resp. right) Dirac
operator of $\Lambda (\R^{3})$:
\begin{equation}
             {\cal D}_{L} \equiv
L({\bf i})\partial _{0}+L({\bf j})\partial _{1}+L({\bf k})\partial _{2} \equiv
                \left( \begin{array}{cccc}
             0   & -\partial _{0}  & -\partial _{1}  & -\partial _{2} \\
  \partial _{0}  & 0              & -\partial _{2}  & \partial _{1} \\
  \partial _{1}  & \partial _{2}  &             0  & -\partial _{0} \\
  \partial _{2}  & -\partial _{1}  & \partial _{0}  &             0
                          \end{array} \right)
\end{equation}
respectively
\begin{equation}
              {\cal D}_{R} \equiv
R({\bf i})\partial _{0}+R({\bf j})\partial _{1}+R({\bf k})\partial _{2} \equiv
               \left( \begin{array}{cccc}
             0   & -\partial _{0}  & -\partial _{1}  & -\partial _{2} \\
  \partial _{0}  & 0              & \partial _{2}  & -\partial _{1} \\
  \partial _{1}  & -\partial _{2}  &             0  & \partial _{0} \\
  \partial _{2}  & \partial _{1}  & -\partial _{0}  &             0
                          \end{array} \right)
\end{equation}
with the properties
$${\cal D}_{L}{\cal D}_{L}^{\ast}=-\bigtriangleup _{3}{\bf 1}_{4} $$
where ${\cal D}^{\ast}=-{\cal D}^{T}$, and respectively
$${\cal D}_{R}{\cal D}_{R}^{\ast}=-\bigtriangleup _{3}{\bf 1}_{4} $$
where: ${\cal D}_{L}^{\ast} = -L({\bf i})\partial _{0}
-L({\bf j})\partial _{1}-L({\bf k})\partial _{2}$ (resp.
${\cal D}_{R}^{\ast} = -R({\bf i})\partial _{0}
-R({\bf j})\partial _{1}-R({\bf k})\partial _{2}$).
Another simple covariant decomposition of the three dimensional Laplacian
$-\bigtriangleup _{3}$ can be described by:
\begin{equation}
  {\cal D}= \left( \begin{array}{cccc}
             0   & \partial _{0}  & \partial _{1}  & \partial _{2} \\
  \partial _{0}  & 0              & -\partial _{2}  & \partial _{1} \\
  \partial _{1}  & \partial _{2}  &             0  & -\partial _{0} \\
  \partial _{2}  & -\partial _{1}  & \partial _{0}  &             0
                          \end{array} \right)
\end{equation}
and
\begin{equation}
{\cal D}^{T}=
               \left( \begin{array}{cccc}
             0   & \partial _{0}  & \partial _{1}  & \partial _{2} \\
  \partial _{0}  & 0              & \partial _{2}  & -\partial _{1} \\
  \partial _{1}  & -\partial _{2}  &             0  & \partial _{0} \\
  \partial _{2}  & \partial _{1}  & -\partial _{0}  &             0
                          \end{array} \right)
\end{equation}
and then ${\cal D}{\cal D}^{T}= \bigtriangleup _{3}{\bf 1}_{4}$.
This corresponds to the choice $a=1,b=1$ and $c=+1$ (resp. $-1$) in
(\ref{allcov}).
The question of the covariance properties of this decomposition was the
starting point of the present research.

\subsection{$D_1\oplus D_1$: Interacting Vector Fields}
The models of this sort describe a doublet of vector fields ${\bf A}=(A_{0},
A_{1},A_{2})$, ${\bf B}=(B_{0},B_{1},B_{2})$ coupled to itself throughout
the noise in the corresponding covariant SPDE.

The realification matrix $E_{(1,1)}$ is chosen to be:
\begin{equation}
E_{(1,1)} =
                    \frac{1}{\sqrt{2}}
              \left( \begin{array}{cccccc}
            i &         0 & -i & 0 &         0 &  0  \\
            1 &         0 &  1 & 0 &         0 &  0  \\
            0 & i\sqrt{2} &  0 & 0 &         0 &  0  \\
            0 &         0 &  0 & i &         0 & -i  \\
            0 &         0 &  0 & 1 &         0 &  1  \\
            0 &         0 &  0 & 0 & i\sqrt{2} &  0
                        \end{array} \right)
\end{equation}
The manifestly real expressions for ${\cal D}_{(1,1)} \in Cov((D_{1}
\oplus D_{1})^{R})$ obtained by the application of $E_{(1,1)}$ are given by:
\begin{equation}
\hat{\cal D} _{(1,1)}=\left( \begin{array}{cccccc}
m_{1} & -aip_{2} & aip_{1} & 0 & -bip_{2} & bip_{1} \\
aip_{2} & m_{1} & -aip_{0} & bip_{2} & 0 & -bip_{0} \\
-aip_{1} & aip_{0} & m_{1} & -bip_{1} & bip_{0} & 0 \\
0 & -cip_{2} & cip_{1} & m_{2} & -dip_{2} & dip_{1} \\
cip_{2} & 0 & -cip_{0} & dip_{2} & m_{2} & -dip_{0} \\
-cip_{1} & cip_{0} & 0 & -dip_{1} & dip_{0} & m_{2}
                        \end{array} \right)
\end{equation}
where the central element $M$ is chosen to $M=m_{1} {\bf 1}_{3}\oplus
m_{2} {\bf 1}_{3}$, $m_{1},m_{2} \in \R$.
\par
Computing $det\hat{\cal D} _{(1,1)}(ip)$ we obtain:
$$
\hbox{det}\hat{\cal D} _{(1,1)}(ip)=m_{1}m_{2}\{ (ad-bc)^{2}p^{4}+
((m_{2}^{2}a^{2})
+2bcm_{1}m_{2}+m_{1}^{2}d^{2})p^{2}+m_{1}^{2}m_{2}^{2} \}.
$$
The conditions for the proper mass spectrum could be easily obtained as
$\hat{{\cal D}}_{(1,1)}(ip)$ is a biquadratic polynom.
Providing that $\hbox{det}\hat{\cal D} _{(1,1)}(ip) \neq 0$ we can invert the
matrix
$\hat{\cal D} _{(1,1)}(ip)$ obtaining the corresponding Green-function:
\begin{equation}
\hat{\cal D}^{-1}_{(1,1)}(p,\alpha)=\frac{1}{\{ f^{2}p^{4}-(h^{2}-2fm_{1}
m_{2})p^{2}+m_{1}^{2}m_{2}^{2} \} }  \left(
                 \begin{array}{c|c}
     G^{(1,1)}_{\mu \nu}(p,\alpha) & G^{(1,2)}_{\mu \nu}(p,\alpha) \\ \hline
     G^{(2,1)}_{\mu \nu}(p,\alpha) & G^{(2,2)}_{\mu \nu}(p,\alpha)
                 \end{array} \right)
\end{equation}
where we have: $f \equiv ad-bc$, $h \equiv am_{2}+dm_{1}$ and
$$G^{(1,1)}_{\mu \nu}(p,\alpha) =
m_{1}m_{2}(-e_{1}p^{2}+m_{1}m_{2}^{2})\delta _{\mu \nu}+
m_{2}(f^{2}p^{2}-m_{2}e_{2})p_{\mu}p_{\nu}+
m_{1}m_{2}(-dfp^{2}+am_{2}^{2})i\varepsilon _{\mu \nu \lambda}p_{\lambda}
$$
$$
G^{(1,2)}_{\mu \nu}(p,\alpha) = bm_{1}m_{2}\{ hp^{2}\delta _{\mu \nu}-
hp_{\mu}p_{\nu}+
(fp^{2}+m_{1}m_{2})i\varepsilon _{\mu \nu \lambda}p_{\lambda} \}
$$
for $\mu ,\nu \in \{ 0,1,2 \}$ with $e_{1} \equiv (d^{2}m_{1}+bcm_{2})$,
$e_{2} \equiv (a^{2}m_{2}+bcm_{1})$. \par
The two last blocks of Green matrix we can obtain by making the following
exchanges: $a \leftrightarrow d$, $m_{1} \leftrightarrow m_{2}$ within the
$G^{(1,1)}_{\mu \nu}(p,\alpha)$ matrix to get $G^{(2,2)}_{\mu \nu}(p,\alpha)$
and
$b \leftrightarrow c$ within $G^{(1,2)}_{\mu \nu}(p,\alpha)$ to get
$G^{(2,1)}_{\mu \nu}(p,\alpha)$. The corresponding two-point Schwinger function
(more precisely the contribution coming from Poisson piece of the noise without
integration over $\nu$ as in the previous case) is given as follows:
\begin{equation}
\hat{S}_{(1,1)}^{(2)}(p,\alpha)=\frac{1}{
   \{ f^{2}p^{4}-(h^{2}-2fm_{1}m_{2})p^{2}+m_{1}^{2}m_{2}^{2} \}^{2}
                 } \left(
                  \begin{array}{c|c}
   \hat{S}^{(1,1)}_{\mu \nu}(p,\alpha) &
                             \hat{S}^{(1,2)}_{\mu \nu}(p,\alpha)  \\ \hline
     \hat{S}^{(2,1)}_{\mu \nu}(p,\alpha) & \hat{S}^{(2,2)}_{\mu \nu}(p,\alpha)
                  \end{array} \right)
\end{equation}
with:
$$
\hat{S}^{(1,1)}_{\mu \nu}(p,\alpha)=
[m_{2}(f^{2}p^{2}-m_{2}e_{2})\alpha_{\lambda}p_{\lambda}-cm_{1}m_{2}h
\beta _{\lambda}p_{\lambda}]^{2}p_{\mu}p_{\nu}+
[m_{2}(f^{2}p^{2}-m_{2}e_{2})\alpha_{\lambda}p_{\lambda}-cm_{1}m_{2}h
\beta _{\lambda}p_{\lambda}]
$$
$$\times
[m_{1}m_{2}(-e_{1}p^{2}+m_{1}m_{2}^{2})(p_{\mu}\alpha _{\nu}+p_{\nu}\alpha
_{\mu})+cm_{1}m_{2}hp^{2}(p_{\mu}\beta_{\nu}+p_{\nu}\beta_{\mu})]+
$$
$$
cm_{1}^{2}m_{2}^{2}hp^{2}(-e_{1}p^{2}+m_{1}m_{2}^{2})(\alpha _{\mu}
\beta_{\nu}+
\alpha _{\nu}\beta _{\mu})+m_{1}^{2}m_{2}^{2}(-e_{1}p^{2}+m_{1}m_{2}^{2})^{2}
\alpha_{\mu}\alpha_{\nu}+(cm_{1}m_{2}hp^{2})^{2}\beta_{\mu}\beta_{\nu},
$$
$$
\hat{S}^{(1,2)}_{\mu \nu}(p,\alpha)=
-m_{1}m_{2}\{ hbm_{2}(f^{2}p^{2}-m_{2}e_{2})(\alpha _{\lambda}p_{\lambda})
^{2}+hcm_{1}(f^{2}p^{2}-e_{1}m_{1})(\beta _{\lambda}p_{\lambda})^{2}-
$$
$$
\alpha _{\lambda}p_{\lambda}\beta _{\omega}p_{\omega}[
(f^{2}p^{2}-m_{2}e_{2})(f^{2}p^{2}-m_{1}e_{1})+m_{1}m_{2}bch^{2}] \}
p_{\mu}p_{\nu}+
m_{1}m_{2}^{2}bhp^{2}[(f^{2}p^{2}-m_{2}e_{2})\alpha _{\lambda}p_{\lambda}-
$$
$$
-m_{1}ch\beta _{\lambda}p_{\lambda}]p_{\mu}\alpha_{\nu}+
m_{1}^{2}m_{2}(-e_{1}p^{2}+m_{1}m_{2}^{2})[(f^{2}p^{2}-m_{1}e_{1})\beta
_{\lambda}p_{\lambda}-m_{2}bh\alpha _{\lambda}p_{\lambda}]p_{\nu}\alpha_{\mu}+
$$
$$
+m_{1}m_{2}^{2}(-e_{2}p^{2}+m_{1}^{2}m_{2})[(f^{2}p^{2}-m_{2}e_{2})
\alpha _{\lambda}p_{\lambda}-m_{1}ch\beta_{\lambda}p_{\lambda}]
p_{\mu}\beta_{\nu}+
m_{1}^{2}m_{2}chp^{2}[(f^{2}p^{2}-m_{1}e_{1})\beta_{\lambda}p_{\lambda}-
$$
$$
-bm_{2}h\alpha_{\lambda}p_{\lambda}]p_{\nu}\beta_{\mu}+
m_{1}^{2}m_{2}^{2}hp^{2}
[b(-e_{1}p^{2}+m_{1}m_{2}^{2})\alpha_{\mu}\alpha_{\nu}+
c(-e_{2}p^{2}+m_{1}^{2}m_{2})\beta _{\mu}\beta _{\nu}]+
$$
$$
+m_{1}^{2}m_{2}^{2}(-e_{1}p^{2}+m_{1}m_{2}^{2})(-e_{2}p^{2}+m_{1}^{2}m_{2})
\alpha_{\mu}\beta_{\nu}+bc(m_{1}m_{2}hp^{2})^{2}\alpha_{\nu}\beta_{\mu}
$$
for $\mu , \nu \in \{ 0,1,2 \}$. \par
We use the notation for the variable of L\'evy measure:
$\alpha \equiv (\alpha _{0}, \alpha _{1} , \alpha _{2},\beta _{0} ,
\beta _{1} , \beta _{2})$.
The block $\hat{S}^{(2,2)}_{\mu \nu}(p,\alpha)$ one can get by the exchanges
$a \leftrightarrow d$, $m_{1} \leftrightarrow m_{2}$ and $b \leftrightarrow c$
in the block $\hat{S}^{(1,1)}_{\mu \nu}(p,\alpha)$ and the block
$\hat{S}^{(2,1)}_{\mu \nu}(p,\alpha)$ by $b \leftrightarrow c$ within
the block $\hat{S}^{(1,2)}_{\mu \nu}(p,\alpha)$.\par
It is worthwile to observe that in the variety of covariant operators
there do not exists a reflection covariant operator.
By specialization of parameters of covariant operator we can find in the
Gaussian part of two-point Schwinger function the Euclidean 2-point function
of two copies of so called Euclidean Proca field introduced in
\cite{Velo,Gross,Yao}.
If we put
\begin{equation}
a=d=0,~ b^{2}=c^{2}=1,~~ bc=-1~~\hbox{and}~~ m_{1}=m_{2}=m
\end{equation}
then we obtain for the Gaussian part
\begin{equation}
\hat{S}_{G;(1,1))}^{(2)}(p)=\left(
                \begin{array}{c|c}
          (\delta _{\mu \nu}+ \frac{p_{\mu}p_{\nu}}{m^{2}})
          \frac{1}{p^{2}+m^{2}} &
                 \begin{array}{ccc}
                       0 & 0 & 0 \\
                       0 & 0 & 0 \\
                       0 & 0 & 0 \\
                 \end{array}       \\ \hline
         \begin{array}{ccc}
                       0 & 0 & 0 \\
                       0 & 0 & 0 \\
                       0 & 0 & 0 \\
                \end{array} & (\delta _{\mu \nu}+ \frac{p_{\mu}p_{\nu}}{m^{2}})
          \frac{1}{p^{2}+m^{2}}
                   \end{array} \right).
\end{equation}
The corresponding covariance matrix $A={\bf 1}_{6}$ (see eq. $2.17$).

\subsection{The $D_{\frac{1}{2}}\oplus D_{\frac{1}{2}}$-case}
This representation seems to be not of physical interest as conflicting the
usual spin-statistic connection. We note that in the case of non positive
quantum field theory the standard spin-statistic theorem could be violated
\cite{Bogoliubov}.
We can use the $D_{\frac{1}{2}}\oplus D_{\frac{1}{2}}$-representation
for the noise $\eta$ transformation rule. In this context the
study of realifications of this representation could be much usefuler than
the analysis of the corresponding covariant operators, Green and Schwinger
functions. However, we mention the case to complete
the list of the lowest dimensional cases.
\par
The chosen realification matrix $E_{(\frac{1}{2},\frac{1}{2})}$
\begin{equation}
E_{(\frac{1}{2},\frac{1}{2})}= \frac{1}{\sqrt{2}}\left(
               \begin{array}{cccc}
               1 & 0 & 0 & 1 \\
               i & 0 & 0 & -i \\
               0 & 1 & -1 & 0 \\
               0 & i & i & 0
               \end{array} \right)
\end{equation}
The covariant operator in the Fourier representation:
$$
   \hat{{\cal D}}_{(\frac{1}{2},\frac{1}{2})}(p)= \left(
                \begin{array}{cc}
 cip_{0}-dip_{1}-aip_{2}+m &-dip_{0}-cip_{1}-bip_{2}  \\
-dip_{0}-cip_{1}+bip_{2} & -cip_{0}+dip_{1}-aip_{2}+m \\
 aip_{0}+bip_{1}+cip_{2} & bip_{0}-aip_{1}-dip_{2}        \\
     -bip_{0}+aip_{1}-dip_{2} & aip_{0}+bip_{1}-cip_{2}
                \end{array} \right.
$$
\begin{equation}
    \left.   \begin{array}{cc}
 aip_{0}-bip_{1}+cip_{2} & bip_{0}+aip_{1}-dip_{2} \\
 -bip_{0}-aip_{1}-dip_{2} & aip_{0}-bip_{1}-cip_{2} \\
-cip_{0}-dip_{1}+aip_{2}+m & dip_{0}-cip_{1}+bip_{2} \\
  dip_{0}-cip_{1}-bip_{2}  & cip_{0}+dip_{1}+aip_{2}+m
                \end{array} \right)
\end{equation}
with $a, b, c, d \in \R $ and $\hbox{det}(\hat{{\cal D}}_{(\frac{1}{2},
\frac{1}{2})}(p))=[(a^{2}+b^{2}+c^{2}+d^{2})p^{2}+m^{2}]^{2}
-4m^{2}b^{2}p^{2}$.
\par
We can get the admissible mass spectrum taking, for example, $b=0$.
\par
We can use the methods presented above to obtain explicit formulas
for the Green functions and the Schwinger functions. The
corresponding expressions are much more complicated than in the examples
before and will therefore not be presented here.\par
\vskip 3 truecm
{\bf Acknowledgments}\par
We thank Professor S.Albeverio for instructive and stimulating
discussions.
We also profited from various discussions with Z.Haba, L.Jak\'obczyk,
and J.L\"offelholz.
The first author thanks the Institute of Theoretical Physics,
University of Wroc\l aw, for the kind hospitality.
The second author thanks Professor Ph.Blanchard for making his stays
at BiBoS Research Centre pleasant and fruitful.
Financial supports from the Bochum-Wroc\l aw exchange programme,
Sonderforschungsbereich 237, the EG Mobility
Grant, and the Polish KBN grant 2PO3B12211
is gratefully acknowledged.

\end{document}